\documentclass[10pt]{article}

\usepackage{amsmath,amssymb,amsthm,amsfonts,amscd,flafter,epsf, epsfig,graphicx,graphs}

\title{Tight contact structures and genus one fibered knots}
\author{John A. Baldwin}
\date{}
\newcommand\x{x}
\newcommand\y{y}
\newcommand\fwa{F_{W_1}}
\newcommand\fwb{F_{W_2}}
\newcommand\fwc{F_{W_3}}
\newcommand\Tight{Tight(\Sigma, \partial \Sigma)}
\newcommand\Dehn{Dehn^+(\Sigma, \partial \Sigma)}
\newcommand\Veer{Veer(\Sigma, \partial \Sigma)}
\newcommand\by{\delta}
\newcommand\w{\x\y^2\x \y^2}
\newcommand\ls{\x^{a_1} \y^{-b_1}...\x^{a_n}\y^{-b_n}}
\newcommand\mcg{\Gamma_\Sigma}
\newcommand\homology{H_1(\Sigma, \mathbb{Z})}
\newcommand\heeg{\widehat{HF}}
\newcommand\wa{\widehat{M}(0;-b_1,..., -b_n)}
\newcommand\ww{\widehat{M}(0;2,2,-b_1,..., -b_n)}
\newcommand\wx{\widehat{M}(-1;2,2,-b_1,..., -b_n)}
\newcommand\wv{\widehat{M}(0;2,2,-b_1,..., -b_n-1)}
\newcommand\Qtwo{Q(b_1+1,b_2,...b_{n-2},b_{n-1}+1)}
\newcommand\Qone{Q(b_1,...,b_n)}
\newcommand\A{A(b_1,...,b_n)}
\newcommand\An{A_{-1}(b_1,...,b_n)}
\newcommand\Az{A_0(b_1,...,b_n)}

\newtheorem{thm}{Theorem}[section]
\newtheorem{lem}[thm]{Lemma}
\newtheorem{claim}[thm]{Claim}
\newtheorem{remark}[thm]{Remark}
\newtheorem{definition}[thm]{Definition}

\begin{document}
\maketitle
\begin{abstract}  
We study contact structures compatible with genus one open book decompositions with one boundary component. Any monodromy for such an open book can be written as a product of Dehn twists around dual non-separating curves in the once-punctured torus. Given such a product, we supply an algorithm to determine whether the corresponding contact structure is tight or overtwisted. We rely on Ozsv{\'a}th-Szab{\'o} Heegaard Floer homology in our construction and, in particular, we completely identify the $L$-spaces with genus one, one boundary component, pseudo-Anosov open book decompositions. Lastly, we reveal a new infinite family of hyperbolic three-manifolds with no co-orientable taut foliations, extending the family discovered in \cite{RSS}. 
 \end{abstract} 

\section{Introduction}
The mapping class group of the once-punctured torus, $\mcg$, is generated by right-handed Dehn twists about dual non-separating curves, $\x$ and $\y$. In an abuse of notation we denote, by $\gamma$, the right-handed Dehn twist around the curve $\gamma \subset \Sigma.$ The left-handed Dehn twist around $\gamma$ is then denoted by $\gamma^{-1}$. When it is unclear whether we are talking about a curve or a Dehn twist, we will use the notation $D_\gamma$ for the right-handed twist around $\gamma$. Given an open-book decomposition ($\Sigma$, $\phi$), we can express $\phi$ as a product of Dehn twists, $\x^{a_1} \y^{b_1}\x^{a_2} \y^{b_2}...\x^{a_n} \y^{b_n}$, with $a_{i}, b_{j} \in \mathbb{Z}$ (In our notation, composition is on the left). There is an equivalence relation on open books given by stabilization/destabilization, and Giroux \cite{Gi}, extending results of Thurston and Winkelnkemper \cite{ThW}, recently showed that equivalence classes of open books are in one-to-one correspondence with isotopy classes of contact structures. Therefore, given a contact structure compatible with an open book, it is natural to ask whether we can infer properties of the contact structure simply by examining its monodromy. For instance, it is well known that if $\phi$ can be expressed as the product of right-handed Dehn twists, then the corresponding contact structure is Stein-fillable \cite{Gi}, hence tight. The converse is also true. In general, however, there are tight contact structures which are not Stein fillable (or even symplectically fillable) \cite{EH}.

Along these lines, we give an algorithm which explicitly determines when a contact structure compatible with a genus one open book with one boundary component is tight and when it is overtwisted. The input to the algorithm is a monodromy, written as a word in the mapping class group in the Dehn twists $\x$ and $\y$. First, we state the main result when $\phi$ is pseudo-Anosov. Recall that pseudo-Anosov is equivalent, in the case of the once-punctured torus, to the condition that $|trace(\phi_{\#})| >2$, where $\phi_{\#}: H_1(\Sigma,\mathbb{Z}) \longrightarrow H_1(\Sigma,\mathbb{Z})$ is the induced map on homology. Pseudo-Anosov monodromies are especially interesting because the mapping torus, $M_{\phi}$, is hyperbolic if $\phi$ is pseudo-Anosov, by a result of Thurston \cite{Th1}. We will return to this fact in section \ref{sec:Lspaces}. From this point forward, $\Sigma$ will denote a genus one surface with one boundary component. The following is merely a preliminary theorem which makes subsequent calculation much easier.

\begin{thm}
\label{thm:WordReduction}
Let $\delta$ be a curve parallel to the boundary, and let $\phi$ be pseudo-Anosov. Then the open book ($\Sigma$, $\phi$) is equal to an open book whose monodromy is of the form 
\begin{description}
\item[I.]$\by^k *\ls$, if $trace(\phi)>2$
\item[II.]$\by^k *\w*\ls$, if $trace(\phi) < -2$
\end{description}
 Here $k \in \mathbb{Z}, a_{i}, b_{j} \in \mathbb{Z}^{\geq 0}$, and $a_i \neq 0 \neq b_j$ for some $i,j$.
\end{thm}

\noindent The pseudo-Anosov version of our main theorem is:

\begin{thm}
\label{thm:TightPseudoAnosov}
Let $\phi$ be a pseudo-Anosov, boundary-fixing automorphism of the once-punctured torus. If $\phi$ is of type I then the contact structure compatible with ($\Sigma$, $\phi$) is tight if and only if $k \geq 1$. Likewise, if $\phi$ is of type II then the contact structure compatible with ($\Sigma$, $\phi$) is tight if and only if $k \geq 0$.
\end{thm}

\noindent We generalize both of these theorems in later sections to account for all monodromies.

To place this result in its proper context, it is necessary to discuss the recent work of Honda, Kazez, and Mati{\'c}. In late 2005, they found a general criterion for the tightness of an open book, introducing the notion of \emph{right-veering} diffeomorphisms \cite{HKM1}. Their result is an improvement over Goodman's sobering arc criterion for overtwistedness \cite{Go}. In particular, they prove that a contact structure $\xi$ is tight if and only if all of its compatible open books $(\Sigma, \phi)$ have right-veering $\phi$. 

In general, however, it is very difficult to prove statements about \emph{all} open books compatible with a given contact structure. Our paper succeeds in characterizing tightness for contact structures in terms of a \emph{single} compatible open book, when the open book has genus equal to one and one boundary component. Honda, Kazez, and Mati{\'c} have independently succeeded in characterizing tightness in terms of a single compatible open book (of genus one with one boundary component). Their results are phrased in terms of the \emph{fractional Dehn twist coefficient} of $\phi$, and can be partially found in \cite{Ho} and \cite{HKM2}. Our approach is entirely different, and it is unique in the fact that we provide an \emph{explicit} algorithm for detecting tightness when given a monodromy written as a product of Dehn twists of the sort described above. Moreover, our method leads very naturally to the discovery of a new family of hyperbolic three-manifolds with no taut foliations. 

Another interesting and related project is to identify those monodromies that give tight contact structures (denote this set by $\Tight$), but which cannot be expressed as the product of right-handed Dehn twists along curves on the once-punctured torus (denote this set by $\Dehn$). This is the central topic of \cite{HKM2}, and it is a step towards characterizing monodromies which produce tight, but non-Stein-fillable contact structures. An advantage of our explicit approach is that it allows us to easily identify a large family of monodromies in $\Tight - \Dehn$. The reader should compare these monodromies with those found in \cite{HKM2}.

\subsection{Organization}
The organization of this paper is as follows: In section \ref{sec:WordReduction}, we prove a generalization of Theorem \ref{thm:WordReduction} by somewhat tedious manipulations in $\mcg$. In section \ref{sec:ContactInvariant}, we calculate the Ozsv{\'a}th-Szab{\'o} contact invariants for the type II monodromies of Theorem \ref{thm:WordReduction} and we prove half of Theorem \ref{thm:TightPseudoAnosov}. In section \ref{sec:Overtwistedness}, we complete the proof of Theorem \ref{thm:TightPseudoAnosov} using Goodman's criterion for overtwistedness. In addition, we generalize Theorem \ref{thm:TightPseudoAnosov}, giving a complete characterization of tightness for all genus one, one boundary component open boooks. In section \ref{sec:TypeD}, we complete the proof of this generalization for type $D$ monodromies, and we discuss $spin^c$ structures. In section \ref{sec:TightDehn}, we analyze $\Tight - \Dehn$. Finally, in section \ref{sec:Lspaces} we classify $L$-spaces which have genus one, one boundary component open book decompositions. This involves a comparison with some of Roberts' results on taut foliations \cite{Ro}, \cite{Ro2}. Moreover, we identify an infinite family of hyperbolic $L$-spaces obtained by surgery on the bindings of these open books. Section \ref{sec:Appendix} is an Appendix containing the proof of Lemma \ref{lem:H1Sum}.

\subsection{Acknowledgements}
I would like to thank Shaffiq Welji, Elisenda Grigsby, Jiajun Wang, Matt Hedden, and Joan Licata for helping me understand Heegaard Floer homology. I am also grateful to John Etnyre and Rachel Roberts for enlightening e-mail discussions. Most of all, I am indebted to my advisor, Peter Ozsv{\'a}th, who suggested that I study the Heegaard Floer homology of these open books. Rarely have I left his office without a new idea.

\newpage
\section{Proof of Theorem 1.1}
\label{sec:WordReduction}
\subsection{Useful Notation and the Mapping Class Group}
As mentioned in the Introduction, the mapping class group, $\mcg$, of the once-punctured torus is generated by right-handed Dehn twists about dual non-separating curves, $\x, \, \y \subset \Sigma$. We orient $x$ and $y$ so that $i(x, y) = +1$, where $i$ is the intersection form on $H_1(\Sigma, \mathbb{Z})$. Given an open book $(\Sigma, \phi)$, where $\phi$ is a word in $\mcg$, it is useful to know how we can change $\phi$ and preserve the open book. We will use the following relations from $\mcg$:

\begin{itemize}
\item{$\x\y\x = \y\x\y$}
\item{$(\w)^2 = (\x\y)^6 = \by$}
\item{If $\gamma$ and $\tau$ are disjoint curves in $\Sigma$, then $\gamma \tau = \tau \gamma$.  }
\end{itemize}
The following notational convention will be useful. Let $M(k; b_1,...,b_n)$ denote the open book $(\Sigma, \by^k * \x \y^{b_1} \x \y^{b_2} ... \x \y^{b_n})$ for any collection of $b_i \in \mathbb{Z}$. Then Theorem \ref{thm:WordReduction} has the following re-formulation: \\

\noindent\textbf{Theorem \ref{thm:WordReduction}.}
\emph{Let $\phi$ be pseudo-Anosov. Then the open book ($\Sigma$, $\phi$) is equal to an open book whose monodromy is of the form }
\begin{description}
\item[I.] $M(k; -b_1,...,-b_n)$, if $trace(\phi)>2$
\item[II.]$M(k; 2,2,-b_1,...,-b_n)$, if $trace(\phi) < -2$
\end{description}
 \emph{Here $k \in \mathbb{Z}, b_{j} \in \mathbb{Z}^{\geq 0}$, and $b_i \neq 0$ for some $i$.} \\

\noindent Below is a list of common "moves" which change the word $\phi$, but preserve the open book $(\Sigma, \phi)$. Each is obtained from a combination of the relations mentioned above, together with the observation that $(\Sigma, w_1 * w_2)$ is the same open book as $(\Sigma, w_2 * w_1),$ where $w_1,w_2$ are words in $\mcg$. This is not a manifestation of relations in $\mcg$, but rather it is due to the fact that the open book $(\Sigma, \phi)$ is constructed from the mapping torus, $M_{\phi}$. 

\begin{lem}
\label{lem:WordMoves}
The following moves preserve the open books:
\begin{enumerate}
\item{$M(k; b_1,...,b_n) = M(k; b_2,...,b_n,b_1)$}
\item{$M(k; b_1,...,\overbrace{b_i,1,b_{i+1}},...,b_n)=M(k; b_1,...,\overbrace{b_i+1,b_{i+1}+1},...,b_n)$}
\item{$M(k; b_1,...,\overbrace{b_i,2,2,2,2,b_{i+1}},...,b_n)=M(k; b_1,...,\overbrace{b_i,1,1,1,1,1,1,b_{i+1}},...,b_n) = M(k+1; b_1,...,b_n) $}
\item{$M(k; b_1,...,\overbrace{b_i,2,b_{i+1}},...,b_n) = M(k; b_1,...,\overbrace{b_i\pm m,2,b_{i+1} \mp m},...,b_n)$}
\item{$M(k; b_1,...,\overbrace{b_i,2,2,b_{i+1}},...,b_n) = M(k; 2,2, b_1,...,\overbrace{b_i,b_{i+1}},...,b_n)$} 
\item{$M(k; b_1,...,\overbrace{b_i,3,b_{i+1}},...,b_n) = M(k; 2,2,b_1,...,\overbrace{b_i-1,b_{i+1}-1},...,b_n)$}
\end{enumerate}
 \end{lem}
 
\noindent \textbf{Proof of Lemma \ref{lem:WordMoves}.}
 $1$, $2$, and $3$ are trivial. $4$ follows because $M(k; b_1,...,b_i,2,b_{i+1},...,b_n) = M(k; b_1,...,b_i-1,1,1,b_{i+1},...,b_n)$ $=$ $M(k; b_1,...,b_i -1,2, b_{i+1}+1,...,b_n)$ on one hand, using moves of type $2$. On the other hand, $M(k; b_1,...,b_i,$ $2,b_{i+1},...,b_n)$ $=$ $M(k; b_1,...,b_i,1,1,b_{i+1}-1,...,b_n)  = M(k; b_1,...,b_i +1,2, b_{i+1}-1,...,b_n)$. $5$ follows from repeated applications of $4$ or, if you prefer, from the fact that $\w$ commutes with $\x$ and $\y$, and thus with everything in $\mcg$. $6$ follows because $M(k; b_1,...,b_i,3,b_{i+1},...,b_n)$ $=$ $M(k; b_1,...,b_i-1,1,2,b_{i+1},...,b_n)$ $=$ $M(k; b_1,...,b_i -1,2,2, b_{i+1}-1,...,b_n)$ $=$ $M(k; 2,2,b_1,...,b_i-1,b_{i+1}-1,...,b_n).$
 \qed

\begin{lem}
\label{lem:WordStandardForm}
Every $(\Sigma,\phi)$ can be expressed as an open book of the form $M(k;b_1,...,b_{2n})$ for $k,b_i \in \mathbb{Z}$. 
\end{lem}

\noindent Note that this is weaker than Theorem \ref{thm:WordReduction} which requires that the $b_i \leq 0$ and some $b_j \neq 0$. \\

\noindent \textbf{Proof of Lemma \ref{lem:WordStandardForm}.}
\begin{claim}
\label{claim:XStandardForm}
 $\x^m = \by^{-1} * \x\y\x\y\x\y\x\y\x\y^{m+1}\x\y$ for $m \in \mathbb{Z}$.
\end{claim}
\noindent Certainly for $m=1$, $$\x = \by^{-1} * \x\y\x\y\x\y\x\y\x\y\x\y\x = \by^{-1} * \x\y\x\y\x\y\x\y\x\y^2\x\y.$$ Now induct: $\x^{m+1} =\x^m*\x$
\begin{eqnarray*}
= &&\by^{-1} * \x\y\x\y\x\y\x\y\x\y^{m+1}\x\y*\by^{-1}*\x\y\x\y\x\y\x\y\x\y^2\x\y \\
= &&\by^{-2}*\x\y\x\y\x\y\x\y\x\y^{m+1}\x\y\x\y\x\y\x\y\x\y\x\y^2\x\y \\
=&&\by^{-2}*\x\y\x\y\x\y\x\y\x\y^{m+1}*\by*\y\x\y\\
=&&\by^{-1} *\x\y\x\y\x\y\x\y\x\y^{m+2}\x\y.
\end{eqnarray*} 
\qed

\noindent Therefore, for any $\phi = \x^{a_1}\y^{b_1}...\x^{a_n}\y^{b_n}$, $(\Sigma,\phi)$ is equal to the open book $M(-n;1,1,1,1,a_1+1,b_1+1,1,1,1,1,a_2+1,b_2+1,...,1,1,1,1,a_n+1,b_n+1)$, completing the proof of Lemma \ref{lem:WordStandardForm}. Note that this can be reduced by the "moves" in Lemma \ref{lem:WordMoves} to $M(\frac{-n}{2};a_1+2,b_1+2,...,a_n+2,b_n+2)$ if $n$ is even, and $M(\frac{-n-1}{2};2,2,a_1+2,b_1+2,...,a_n+2,b_n+2)$ if $n$ is odd.

\subsection{Proof of Theorem \ref{thm:WordReduction}.}

\begin{remark}
\label{remark:PartialWordReduction}
We can use Claim \ref{claim:XStandardForm} to give another "move" among open books:
$$M(k;b_1,...,b_{2n}) = M(k-2n;1,1,1,1,2,b_1+1,...,1,1,1,1,2,b_{2n}+1) = M(k-n;3,b_1+2,...,3,b_{2n}+2).$$

\noindent Any open book can be described by $M(k;b_1,...b_{2n})$ by Lemma \ref{lem:WordStandardForm}, and repetitions of the above move show that $M(k;b_1,...b_{2n}) = M(d;p_1,p_2,...,p_m)$ for some $d$ and some collection of $p_i \geq 3$.  
\end{remark}

\begin{lem}
\label{lem:PartialWordReduction}
$M(k;b_1,...,b_n)$ can be written as one of the following types of open books, via the word moves detailed above:

\begin{description}
\item[1.]$M(d;p_1,...,p_m)$, where the $p_i \geq 4$
\item[2.]$M(d; 2,2, p_1,...,p_m)$, where the $p_i \geq 4$
\item[3.]$M(d;2,p_1)$, where $p_1 \geq 2$
\item[4.]$M(d;2,2,2,p_1)$, where $p_1 \geq 2$
\item[5.]$M(d;p_1)$, where $p_1 \geq1$
\item[6.]$M(d;2,2,p_1)$, where $p_1 \geq 1$
\end{description}
\end{lem}

\noindent \textbf{Proof of Lemma \ref{lem:PartialWordReduction}.}
Suppose that our open book is of the form $M(d;p_1,...,p_m)$ (or $M(d;2,2,p_1,...,p_m)$) where the $p_i \geq 3$. Any open book can be written in this form by Remark \ref{remark:PartialWordReduction}. We perform the following iteration:

\begin{description}
\item[Step 1:] If all of the $p_i \geq 4$, then we can stop as we are in type $1$ (or $2$). If one of the $p_i = 3$, and $m \geq 3$, then go to Step 2. Otherwise, stop.
\item[Step 2:] Apply the word move in Lemma \ref{lem:WordMoves}(6) to obtain $M(d;2,2,p_1,...,p_{i-1}-1,p_{i+1}-1,...,p_m)$ (or $M(d+1;p_1,...,p_{i-1}-1,p_{i+1}-1,...,p_m))$.  If neither $p_{i-1}-1$ nor $p_{i+1}-1$ is equal to $2$ then iterate, i.e., go to Step 1. If both $p_{i-1}-1,p_{i+1}-1 = 2$, and $m \geq 4$ then go to Step 3. If exactly one of $p_{i-1}-1,p_{i+1}-1$ is $2$, say $p_{i-1} -1= 2$, and $m \geq 4$, then go to Step 4. Otherwise, stop.
\item[Step 3:] Apply the word moves in Lemma \ref{lem:WordMoves} to obtain $M(d+1;p_1,...,p_{i-2},p_{i+2},...,p_m)$ (or $M(d+1;2,2,p_1,...,p_{i-2},p_{i+2},...,p_m))$. Return to Step 1.
\item[Step 4:] Our open book is of the form $M(d;2,2,p_1,...,2,p_{i+1}-1,...,p_m)$ (or $M(d+1;p_1,...,2,p_{i+1}-1,...,p_m))$ = $M(d;2,2,p_1,...,2,2,p_{i+1}-1+p_{i-2}-2,...,p_m)$ = $M(d+1;p_1,...,p_{i+1}-1+p_{i-2}-2,...,p_m)$ (or $M(d+1;2,2,p_1,...,p_{i+1}-1+p_{i-2}-2,...,p_m)).$ Return to Step 1.
\end{description}

All that remains of the proof is to show that this iteration stops exactly when we are in one of the types of Lemma \ref{lem:PartialWordReduction}. It stops at Step 1 if we ever reach a point in which our open book is $M(d;p_1,...,p_m)$ (or $M(d;2,2,p_1,...,p_m)$) where the $p_i \geq 4$. Here, we are in type $1$ (or $2$). It stops at Step 1 if one of the $p_i =3$, and $m = 2$. In this case, the open book looks like $M(d;3,p_2)$ (or $M(d;2,2,3,p_2))$  = $M(d;2,1,p_2-1) = M(d;p_2-1,2,1) = M(d;p_2-2,2,2) = M(d;2,2,p_2-2)$ (or $M(d+1;p_2-2)$), and we are in type $6$ (or $5$).  It stops at Step 1 if one of the $p_i =3$, and $m = 1$, in which case we are in type 5 or 6. It stops at Step 2 if $p_i = 3$, both $p_{i-1}-1, p_{i+1}-1 = 2$, and $m=3$. In this case, the open book looks like $M(d;2,2)$ (or $M(d;2,2,2,2)$), and we are in type $3$ (or $4$). Finally, it stops at Step 2 if one of the $p_i = 3$, exactly one of $p_{i-1}-1, p_{i+1}-1$ is equal to $2$ (say $p_{i-1}-1 = 2$), and $m = 3$. In this case, the open book looks like $M(d;2,p_{i+1}-1)$ (or $M(d;2,2,2,p_{i+1}-1)$), and we are in type $3$ (or $4$). 
\qed \\

The following theorem generalizes Theorem \ref{thm:WordReduction}.

\begin {thm}
\label{thm:WordReductionGeneral}
Any open book can be written as $(\Sigma, \phi)$, where $\phi$ is one of the following types:
\begin{description}
\item[A.] $\by^d *\x^{a_1}\y^{-1}...\x^{a_n}\y^{-1}$, where the $a_i \geq 0$, some $a_j \neq 0$.
\item[B.] $\by^d * \w*\x^{a_1}\y^{-1}...\x^{a_n}\y^{-1}$, where the $a_i \geq 0$, some $a_j \neq 0$.
\item[C.] $\by^d *\y^m$, for $m \in \mathbb{Z}$.
\item[D.] $\by^d *\w*\y^m$, for $m \in \mathbb{Z}$.
\item[E.] $\by^d *\x^m\y^{-1}$, where $m \in \{-1,-2,-3\}$
\item[F.] $\by^d *\w*\x^m\y^{-1}$, where $m \in \{-1,-2,-3\}$
\end{description}
Only types A and B are pseudo-Anosov.
\end{thm}

\noindent \textbf{Proof of Theorem \ref{thm:WordReductionGeneral}.}
We simply need to show that monodromies of the open books in types 1-6 of Lemma \ref{lem:PartialWordReduction} can be expressed as monodromies of types A-F. From Claim \ref{claim:XStandardForm}, for $m$ even, we can write 

$$\x\y^{p_1}...\x\y^{p_m} = \by^{\frac{m}{2}}* \y^{-1}\x^{p_1-3}\y^{-1}\x^{p_2-4}...\y^{-1}\x^{p_{m-1}-4}\y^{-1}\x^{p_m-3}\y^{-1}\x^{-1}$$
\noindent And for $m$ odd, we can write

$$\x\y^{p_1}...\x\y^{p_m} = \by^{\frac{m+1}{2}}* \y^{-1}\x^{p_1-3}\y^{-1}\x^{p_2-4}...\y^{-1}\x^{p_{m-1}-4}\y^{-1}\x^{p_{m}-3}\y^{-1}\x^{-1}\y^{-1}\x^{-1}\y^{-1}\x^{-1}\y^{-1}\x^{-1}.$$

\noindent Then we can substitute these identities into to the monodromies of types 1,2,5,6, perform the necessary word moves, and see which of the types A-F we get. For monodromies of types 3 and 4, the reduction is easier. Observe that $\x\y^2\x\y^{p_1} =  \y^0\x\y^2\x\y^{p_1} =  \y^2\x\y^2\x*\y^{p_1-2} = \x\y^2\x\y^2*\y^{p_1-2}$, and substitute. We give the results of these substitutions below.

\begin{description}
\item[1.]$M(d;p_1,...,p_m) = \by^{d+\frac{m}{2}}*\x^{p_1-4}\y^{-1}...\x^{p_m-4}\y^{-1}$ if $m$ is even; $\by^{d+\frac{m-1}{2}}*\w*\x^{p_1-4}\y^{-1}...\x^{p_m-4}\y^{-1}$ if $m$ is odd.
\item[2.]$M(d; 2,2, p_1,...,p_m) = \by^{d+\frac{m}{2}}*\w*\x^{p_1-4}\y^{-1}...\x^{p_m-4}\y^{-1}$ if $m$ is even; $\by^{d+\frac{m+1}{2}}*\x^{p_1-4}\y^{-1}...\x^{p_m-4}\y^{-1}$ if $m$ is odd.
\item[3.]$M(d;2,p_1) = \by^d*\w*y^{p_1-2}$.
\item[4.]$M(d;2,2,2,p_1) = \by^{d+1}*y^{p_1-2}$. 
\item[5.]$M(d;p_1) = \by^d*\w*x^{p_1-4}\y^{-1}$. 
\item[6.]$M(d;2,2,p_1) = \by^{d+1}*\x^{p_1-4}\y^{-1}$.
\end{description}

\qed \\

\noindent \textbf{Proof of Theorem \ref{thm:WordReduction}.}
The proof of Theorem \ref{thm:WordReduction} is finished by the observation that only monodromies of types $A$ and $B$ are pseudo-Anosov. This is seen by computing $trace(\phi_{\#})$ for each of these types. If we let $([x],[y])$ be our basis for $\homology$, then 

$$(\x^{m})_{\#} = 
\left[
\begin {array}{cc}
1&m\\
0&1
\end {array}
\right] \;\;\; (\y^{m})_{\#} = 
\left[
\begin {array}{cc}
1&0\\
-m&1
\end {array}
\right] \; \;\;(\by^{m})_{\#} = 
\left[
\begin {array}{cc}
1&0\\
0&1
\end {array}
\right] \;\;\; (\w)_{\#} = 
\left[
\begin {array}{cc}
-1&0\\
0&-1
\end {array}
\right]$$ 

\noindent The various word moves that we have illustrated certainly preserve $trace(\phi_{\#})$. So, for an arbitrary $\phi \in \mcg$ it is clear that $trace(\phi_{\#})>2$ if and only if $(\Sigma, \phi)$ is equal to an open book with type $A$ monodromy. Likewise, $trace(\phi_{\#})<-2$ if and only if $(\Sigma, \phi)$ is equal to an open book with type $B$ monodromy. Since these type $A$ and $B$ monodromies have the properties required by Theorem \ref{thm:WordReduction}, the proof is complete. 
\qed

\begin{remark}
\label{remark:Classification}
Note that types $C$ and $D$ are reducible, while types $E$ and $F$ are periodic.
\end{remark}

\newpage
\section{Computing the Contact Invariants}
\label{sec:ContactInvariant}
\subsection{The Contact Invariant and Surgery Exact Triangles}
To any contact three-manifold $(Y,\xi)$, we can associate a class, $c(\xi) \in \widehat{HF}(-Y)/\pm{1}$, which is an invariant of the contact structure $\xi$ up to isotopy \cite{OSz1}. We will be using $\mathbb{Z}_2$ coefficients throughout to avoid ambiguity in sign. This invariant encodes information related to the tightness of $\xi$. For instance, Ozsv{\'a}th and Szab{\'o} prove that if $\xi$ is overtwisted, then $c(\xi) = 0$. On the other hand, if $\xi$ is Stein fillable, then $c(\xi) \neq 0$. The precise relationship between $c(\xi)$ and the tightness of $\xi$ is still unknown - there are tight contact structures with vanishing contact invariant \cite{Gh}. Yet, we show, in the case of contact structures compatible with genus one, one boundary component open books, that $c(\xi) = 0$ if and only if $\xi$ is overtwisted. The contact invariant is defined in terms of a compatible open book decomposition, $c(\Sigma, \phi)$, and it satisfies the following property \cite{OSz1}: 

\begin{thm} 
\label{thm:NaturalityOfContactInvariant}
If $(\Sigma, \phi)$ is an open book decomposition for Y, and $\gamma \subset Y - B$ is a curve supported in a page of the open book (B is the binding), which is not homotopic to the boundary, then $(\Sigma,  t_{\gamma}^{-1}*\phi )$ induces an open book decomposition of $Y_{+1}(\gamma)$, and under the map

$$F_{W}: \heeg(-Y) \longrightarrow \heeg(-Y_{+1}(\gamma))$$

\noindent obtained by the two-handle addition (and summing over all $spin^{c}$ structures), we have that 

$$F_{W}(c(\Sigma, \phi)) = \pm c(\Sigma, t_{\gamma}^{-1} * \phi).$$

\end{thm}

\noindent In particular, this tells us that if $c(\Sigma, \phi) \neq 0$, then $c(\Sigma, t_{\gamma}*\phi) \neq 0$, where $(\Sigma,  t_{\gamma}*\phi)$ is the induced open book decomposition of $Y_{-1}(\gamma)$.

In this section, we compute the contact invariants $c(\Sigma,\phi)$, where $\phi$ is a monodromy of type II in Theorem \ref{thm:WordReduction}. Borrowing the notation from section \ref{sec:WordReduction}, we have:

\begin{thm}
\label{thm:ContactInvariant}
If $b_{j} \in \mathbb{Z}^{\geq 0}$ and some $b_i \neq 0$ then $c(M(0;2,2,-b_1,...,-b_n)) \neq 0$.
\end{thm}

\noindent Recall that in $\mcg$, $\by = (\x\y)^6 = (\w)^2$. Then Theorems \ref{thm:ContactInvariant} and  \ref{thm:NaturalityOfContactInvariant} imply that for $k \geq 1$ $$c(M(k;-b_1,...,-b_n)) \neq 0$$
\noindent Hence, the contact structure compatible with $M(k;-b_1,...,-b_n)$ is tight for $k \geq 1$. In exactly the same way $M(k;2,2,-b_1,...,-b_n)$ is tight for $k \geq 0$. This proves half of Theorem \ref{thm:TightPseudoAnosov}.

Before we give the proof of Theorem \ref{thm:ContactInvariant}, we must examine the Heegaard Floer homology of the three-manifolds underlying these open books. In the language of Heegaard Floer homology, a rational homology 3-sphere $Y$ is an $L$-space if $\widehat{HF}(Y,s) \cong \mathbb{Z}$ for every $spin^{c}$-structure $s$ on $Y$. $L$-spaces are closed under surgeries in the following sense \cite{OSz4}: 

\begin{thm} 
\label{thm:LspaceH1}
Suppose $K \subset Y$ is a knot with framing $\lambda$ and 

$$|H_1(Y_{\lambda}(K))| = |H_1(Y)| + |H_1(Y_{\mu+\lambda}(K))|$$

\noindent where $\mu$ denotes the meridian for the knot, and $Y_{\lambda}(K)$ is the 3-manifold obtained from $Y$ by performing surgery on $Y$ along $K$ with framing $\lambda$. Then, if $Y$ and $Y_{\mu+\lambda}(K)$ are both $L$-spaces, so is $Y_{\lambda}(K)$. Furthermore, the map 

$$F_{W}: \heeg(Y) \longrightarrow \widehat{HF}(Y_{\lambda}(K))$$

\noindent obtained by the two-handle addition (and summing over all $spin^c$ structures) is injective.
\end{thm}

\begin{figure}[htbp]
\begin{center}
\includegraphics[height=9cm]{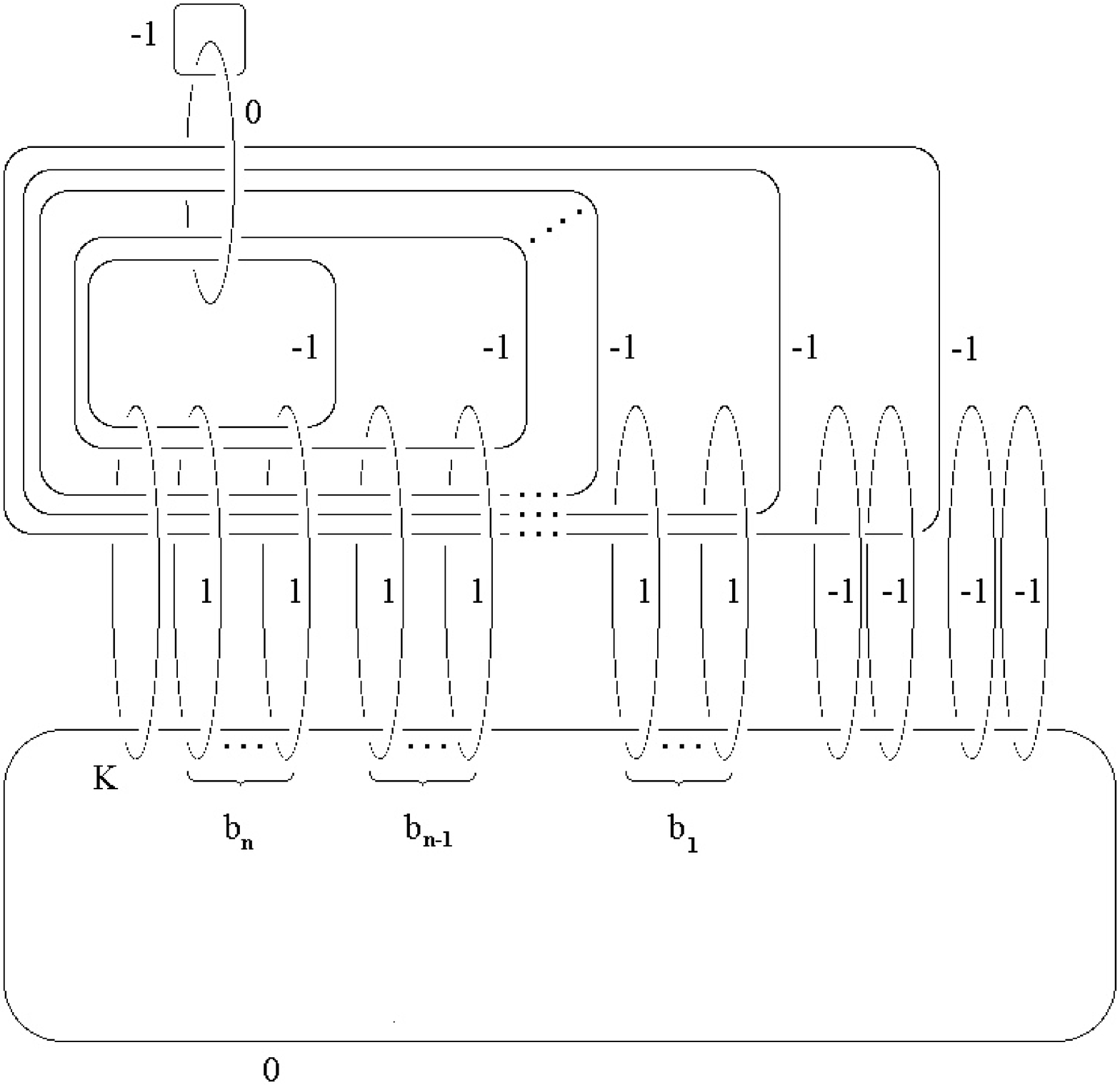}
\caption{\quad Surgery diagram for the 3--manifold $\ww$}
\label{fig:Kirbyww}
\end{center}
\end{figure}

Let $\widehat{M}(k;a_1,...,a_m)$ denote the oriented three-manifold underlying the open book $M(k;a_1,...,a_m)$. Then $\ww$ has the surgery diagram illustrated in Figure \ref{fig:Kirbyww}. The knot $K$ in Figure \ref{fig:Kirbyww} lies in a page of the open book $M(0;2,2,-b_1,...,-b_n)$, and $+1$ surgery on $K$ yields the manifold $\wv$ with induced open book decomposition $M(0;2,2,-b_1,...,-b_n-1).$ Then, by Theorem \ref{thm:ContactInvariant}, the map

$$F_{W}: \heeg(-\ww) \longrightarrow \heeg(-\wv)$$

\noindent takes $$c(M(0;2,2,-b_1,...,-b_n)) \mapsto c(M(0;2,2,-b_1,...,-b_n-1)).$$ The manifolds $-\widehat{M}(0;2,2,-b_1,...,-b_n)$ and $-\widehat{M}(0;2,2,-b_1,...,-b_n-1)$ fit into a surgery exact triangle: \\

\[
\begin{graph}(6,2)
\graphlinecolour{1}\grapharrowtype{2}
\textnode {A}(-1,1.5){$\heeg(-\ww)$}
\textnode {B}(7, 1.5){$\heeg(-\wv) $}
\textnode {C}(3, -.5){$\heeg(\Qtwo)$}
\diredge {A}{B}[\graphlinecolour{0}]
\diredge {B}{C}[\graphlinecolour{0}]
\diredge {C}{A}[\graphlinecolour{0}]
\freetext (3,1.8){$F_W$}
\end{graph}
\]
\\

\begin{figure}[!htbp]
\begin{center}
\includegraphics[height=8cm]{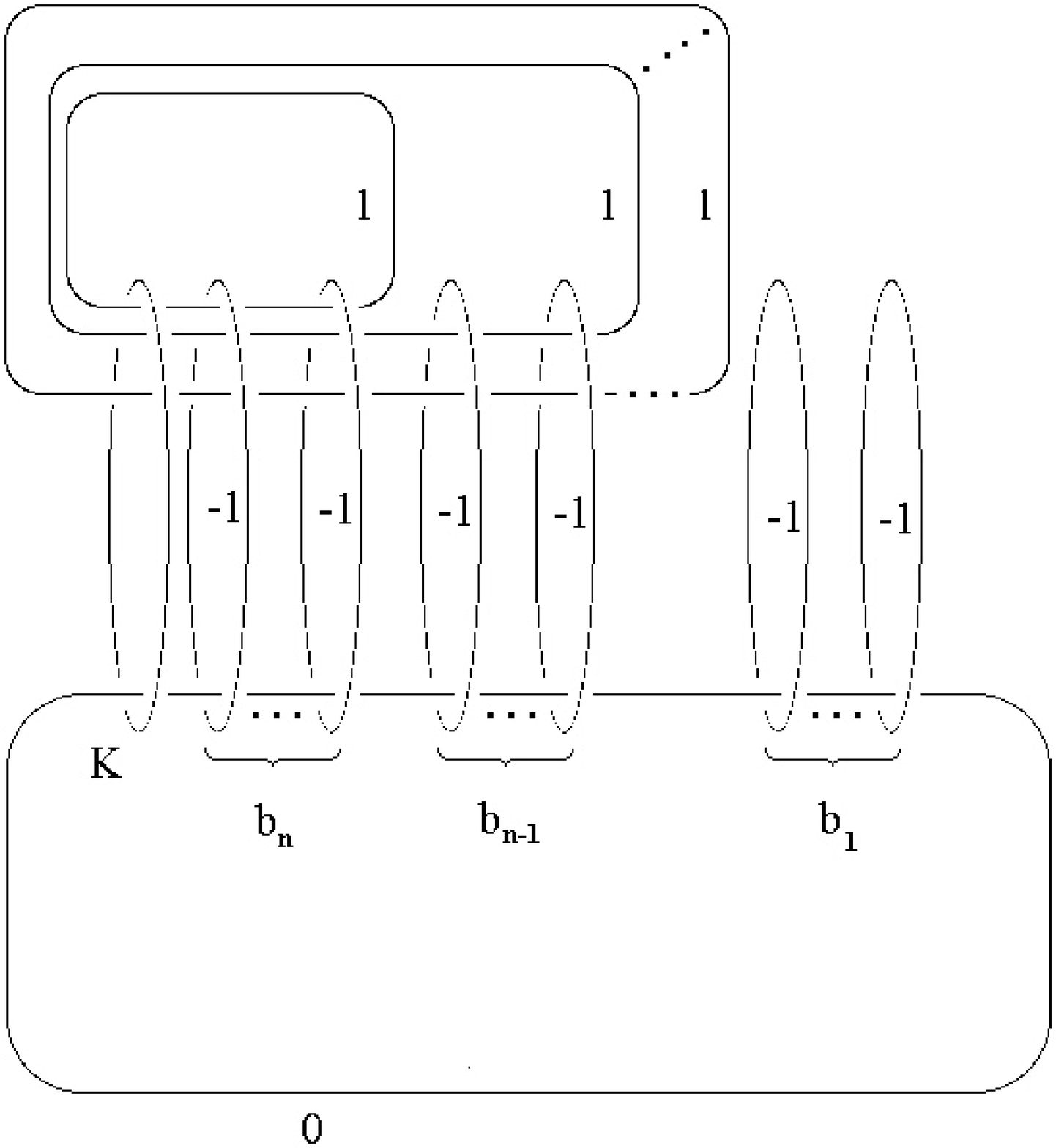}
\caption{\quad Surgery diagram for the 3--manifold $\Qone$}
\label{fig:KirbyQone}
\end{center}
\end{figure}

\noindent $\Qone$ is defined to be the three-manifold given by the surgery diagram in Figure \ref{fig:KirbyQone}. After a sequence of blowdowns and handleslides it can be shown that $\Qtwo$ is the manifold obtained by 0-surgery on the knot $K$ in $-\wv$. 


Our calculation of the contact invariants depends on the following three lemmas:

\begin{lem}
\label{lem:QLspace}
For $b_{j} \in \mathbb{Z}^{\geq 0}$ and some $b_i \neq 0$, $\Qone$ is an L-space.
\end{lem}

\begin{lem}
\label{lem:wwLspace}
For $b_{j} \in \mathbb{Z}^{\geq 0}$ and some $b_i \neq 0$, $-\ww$ is an L-space.
\end{lem}

\begin{lem}
\label{lem:H1Sum}
For $b_{j} \in \mathbb{Z}^{\geq 0}$ and some $b_i \neq 0$, 
$$|H_1(-\wv)|= |H_1(-\ww)|+|H_1(\Qtwo)|.$$
\end{lem}

\noindent Theorem \ref{thm:ContactInvariant} follows immediately: \\

\noindent \textbf{Proof of Theorem \ref{thm:ContactInvariant}.}
$M(0;2,2,-1,\overbrace{0,...,0}^{m}) = M(0;2,2,0-1,0,...,0) = ... = M(0;2,2,0,...,0,-1)$ is the open book given by the monodromy $\w*\x\y^{-1}\x^m = \x\x\y\x\y\x\y*\y^{-1}\x^m = \x^2\y\x\y\x^{m+1} = \x^3\y\x^{m+2}$. And $(\Sigma, \x^3\y\x^{m+2}) = (\Sigma, \y\x^{m+5})$. But then $M(0;2,2,-1,0,...,0)$ is Stein-fillable since it can be written as an open book whose monodromy consists solely of right-handed Dehn twists. Hence, $c(M(0;2,2,-1,0,...,0)) \neq 0.$ Lemmas \ref{lem:QLspace} - \ref{lem:H1Sum}, together with Theorem \ref{thm:LspaceH1} imply that $F_W$ is injective. Consequently $$c(M(0;2,2,-b_1,..., -b_n)) \neq 0 \Longrightarrow c(M(0;2,2,-b_1,..., -b_n-1)) \neq 0.$$ Since $c(M(0;2,2,-1,0,...,0)) \neq 0$, we can induct (as in the later proofs of Lemmas \ref{lem:QLspace} and \ref{lem:wwLspace}) to show that $c(M(0;2,2,-b_1,..., -b_n))\neq 0$ for $b_{j} \in \mathbb{Z}^{\geq 0}$ and some $b_i \neq 0$, completing the proof of Theorem \ref{thm:ContactInvariant}.
\qed

\subsection{Proof of Lemmas \ref{lem:QLspace} - \ref{lem:H1Sum}.}
Lemma \ref{lem:H1Sum} follows from direct computation: compare the determinants of the linking matrices of the surgery diagrams for the manifolds $-\ww$, $-\wv$, and $\Qtwo$. We save this proof for the Appendix.  \\

\noindent \textbf{Proof of Lemma \ref{lem:wwLspace}.}
In the proof of Theorem \ref{thm:ContactInvariant}, we noted that $M(0;2,2,-1,\overbrace{0,...,0}^{m})$ is the open book $(\Sigma, \y\x^{m+5})$ for some $m$. But then, $$-\widehat{M}(0;2,2,-1,0,...,0) = -\widehat{M}(0;2,2,0-1,0,...,0) = ... = -\widehat{M}(0;2,2,0,...,0,-1) = -L(m+5,1)$$ which is an $L$-space. Our proof proceeds by induction on the $b_i$. Suppose that $b_i \geq 1$ for some $i$ and either $b_j >1$ for some $j$, or $b_k \neq 0$ for some $k \neq i$ (otherwise, $-\ww = -\widehat{M}(0;2,2,-1,0,...,0)$), then $$-\ww = -\widehat{M}(0;2,2,-b_{i+1}, ...,-b_{n}, -b_1,...,-b_i)$$ and by induction, we know that $-\widehat{M}(0;2,2,-b_{i+1}, ...,-b_{n}, -b_1,...,-(b_i-1))$ is an $L$-space. Since $Q(b_{i+1}+1, ...,b_{n}, b_1,...,b_{i-1}+1)$ is an $L$-space we can conclude, by Lemma \ref{lem:QLspace} and Theorem \ref{thm:LspaceH1}, that $-\ww = -\widehat{M}(0;2,2,-b_{i+1}, ...,-b_n, -b_1,...,-b_i)$ is an $L$-space.
\qed \\

\noindent \textbf{Proof of Lemma \ref{lem:QLspace}.}
This proof is virtually identical in technique to that of Lemma \ref{lem:wwLspace}. $Q(b_1,...,b_{n-1},b_n+1)$ is the manifold obtained by performing $-1$-surgery on the knot $K$ in $\Qone$ (see Figure \ref{fig:KirbyQone}). Meanwhile, $Q(b_1,...,b_{n-2},b_{n-1}+1)$ is the manifold obtained by performing $0$-surgery on the knot $K$. The following claim is the analogue of Lemma \ref{lem:H1Sum} for these manifolds $Q$. 

\begin{claim}
\label{claim:QH1Sum}
For $b_{j} \in \mathbb{Z}^{\geq 0}$ and some $b_i \neq 0$ then
$$|H_1(Q(b_1,...,b_{n-1},b_n+1))| = |H_1(\Qone)|+|H_1(Q(b_1,...,b_{n-2},b_{n-1}+1))|.$$
\end{claim}

\noindent Again, this is proved directly by comparing the determinants of the linking matrices for these three manifolds. By Theorem \ref{thm:LspaceH1}, if $Q(b_1,...,b_{n-2},b_{n-1}+1)$ and $\Qone$ are $L$-spaces, then so is $Q(b_1,...,b_{n-1},b_n+1)$. To complete the proof of Lemma \ref{lem:QLspace}, we proceed by induction, as before. 

Observe that $Q(0,a_1,...,a_m) = Q(a_1,...,a_m) = Q(a_1,...,a_m,0)$, and $Q(1) = S^3$, an $L$-space. Now, we induct on $n$ and on $b_n$ simultaneously. Consider $\Qone$. If $b_n = 0$ (or $b_1 = 0$), then $\Qone = Q(b_1,...,b_{n-1})$ (or $Q(b_2,...,b_n)$) and we can conclude  by our induction on n that $\Qone$ is an $L$-space. If $b_n \neq 0$ and $b_1 \neq 0$, then by our induction on $b_n$, $Q(b_1,...,b_{n-1},b_n-1)$ is an $L$-space; and by induction on $n$, $Q(b_1,...,b_{n-2},b_{n-1}+1)$ is an $L$-space. Combining this with Claim \ref{claim:QH1Sum} and Theorem \ref{thm:LspaceH1}, we can conclude that $\Qone$ is also an $L$-space. This completes the proof of Lemma \ref{lem:QLspace}, and, consequently, ties up the remaining loose end in the proof of Theorem \ref{thm:ContactInvariant} as well as half of Theorem \ref{thm:TightPseudoAnosov}.
\qed

\begin{remark}
\label{remark:OtherLspaces}
Using this technique, we can prove that the manifolds $\widehat{M}(0;-b_1,...,-b_n)$ and $\ww$ are also $L$-spaces for $b_{j} \in \mathbb{Z}^{\geq 0}$ and some $b_i \neq 0$.
\end{remark}

\newpage
\section{Overtwistedness and Sobering Arcs}
\label{sec:Overtwistedness}
In this section, we prove the second half of Theorem \ref{thm:TightPseudoAnosov}; that is, for $b_{j} \in \mathbb{Z}^{\geq 0}$ and some $b_i \neq 0$, $M(k;-b_1,...,-b_n)$ is overtwisted if $k < 1$ and $M(k;2,2,-b_1,...,-b_n)$ is overtwisted if $k < 0$. These statements follow directly from Goodman's sobering arc criterion for overtwistedness \cite{Go}. First, a bit of background material. 

Given two properly embedded oriented arcs $\alpha$, $\beta$ with common boundary points in $\Sigma$, let $\beta'$ be an arc transverse to $\alpha$ that minimizes intersections with $\alpha$ over boundary-fixing isotopies of $\beta$. Then let $i_{alg}(\alpha, \beta)$ denote the oriented intersection number of $\alpha$ with $\beta'$, summed over points in the interiors of the arcs. Let $i_{geom}(\alpha,\beta)$ be the unsigned number of interior intersection points of $\alpha$ and $\beta'$. And let $i_\delta (\alpha,\beta)$ be one-half the oriented sum of intersections at the boundaries of the arcs $\alpha$ and $\beta'$. In our case, suppose $\alpha$ is an arc in the page $\Sigma$ of an open book, $(\Sigma,\phi)$. Give the arc $\phi(\alpha)$ the orientation which is opposite the pushed-forward orientation of $\alpha$. 

\begin{definition}
\label{def:SoberingArc}
A properly embedded arc $\alpha \subset \Sigma$ is sobering for the monodromy $\phi$ if $$i_{alg}(\alpha, \phi(\alpha)) + i_{geom}(\alpha,\phi(\alpha)) + i_\delta (\alpha,\phi(\alpha))\leq 0$$ and $\alpha$ is not isotopic to $\phi(\alpha)$.
\end{definition}

\begin{figure}[!htbp]
\begin{center}
\includegraphics[height=8cm]{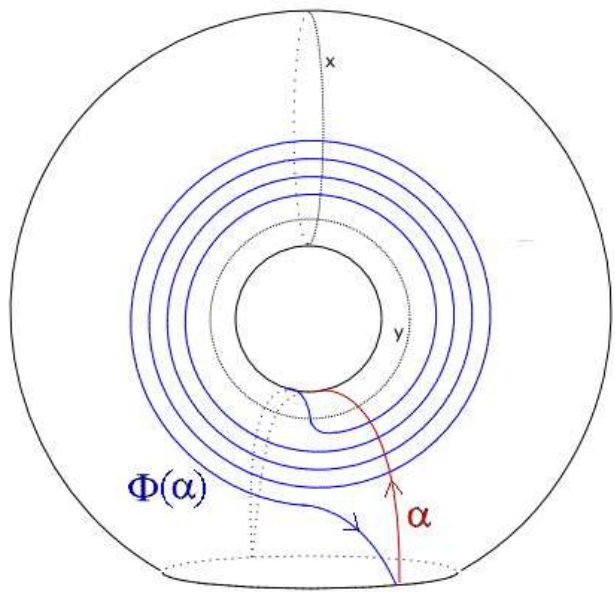}
\caption{\quad A sobering arc in $\alpha \subset \Sigma$ for the monodromy $\phi = \y^{-4}.$ Here $y$ is a longitude and $x$ a meridian.}
\label{fig:SoberingArc}
\end{center}
\end{figure}

\noindent See Figure \ref{fig:SoberingArc} for an illustration of sobering arcs. Goodman showed that

\begin{thm}
\label{thm:OvertwistedSoberingArc}
If there is a sobering arc $\alpha \subset \Sigma$ for the monodromy $\phi$, then the open book $(\Sigma, \phi)$ is overtwisted.
\end{thm}

\noindent \textbf{Proof of Theorem \ref{thm:TightPseudoAnosov}.}
For $k <1$, the open book $M(k;-b_1,...,-b_n)$ has monodromy which is the composition of left-handed Dehn twists around the curve $y$ with Dehn twists around the curve $x$. Therefore, it is clear that an arc $\alpha$ which crosses the curve $y$ once is sobering, as $i_{alg}(\alpha, \phi(\alpha)) + i_{geom}(\alpha,\phi(\alpha)) = 0$, and $i_\delta (\alpha,\phi(\alpha)) = -1.$ In general, whenever a monodromy consists of left-handed Dehn twists around $y$ and arbitrary Dehn twists around $x$ (or vice versa), there is an obvious sobering arc. This is the case for the monodromy of the open book $M(k;2,2,-b_1,...,-b_n)$ if $k < 0$. \qed \\

We generalized Theorem \ref{thm:WordReduction} in Theorem \ref{thm:WordReductionGeneral}. Here we give the corresponding generalization of Theorem \ref{thm:TightPseudoAnosov}.

\begin{thm}
\label{thm:TightGeneral}
The following is a complete classification of tight contact structures compatible with genus one, one boundary component open books. A-F correspond to the monodromies of Theorem \ref{thm:WordReductionGeneral}.
\begin{description}
\item[A:] Tight if and only if $d \geq 1$.
\item[B:] Tight if and only if $d \geq 0$.
\item[C:] Tight if and only if either $d > 0$ or $d = 0$ and $m \geq 0$. 
\item[D:] Tight if and only if $d \geq 0$.
\item[E:] Tight if and only if $d \geq 1$.
\item[F:] Tight if and only if $d \geq 0$. 
\end{description}
\end{thm}

Theorem \ref{thm:TightGeneral} is a complete classification of tight contact structures compatible with genus one, one boundary component open books. \\

 \noindent \textbf{Proof of Theorem \ref{thm:TightGeneral}.}
 For types $A$ and $B$, this statement is simply Theorem \ref{thm:TightPseudoAnosov}. The proof for types $C-F$ is a combination of Goodman's sobering arc criterion for overtwistedness (Theorem \ref{thm:OvertwistedSoberingArc}), and Theorems \ref{thm:ContactInvariant}, \ref{thm:NaturalityOfContactInvariant}, and \ref{thm:TightPseudoAnosov}.
\begin{description}

\item[C:] If $d>0$, then $\by^d*\y^m = \by^{d-1}*\y^2\x\y^2*\x\y^2\x\y^2*\x\y^m$. Since $c(\Sigma, \x\y^2\x\y^2*\x\y^m) \neq 0$ for any $m$ (by Theorem \ref{thm:ContactInvariant}), we conclude by Theorem \ref{thm:NaturalityOfContactInvariant} that $c(\Sigma, \by^d*\y^m) \neq 0$, hence this open book is tight. \\

If $d < 0$, then $\by^d*\y^m$ consists of left handed Dehn twists around $x$ with Dehn twists around $y$, and is therefore overtwisted by Theorem \ref{thm:OvertwistedSoberingArc}. \\

If $d = 0$, then $c(\Sigma, \y^m) \neq 0$ if $m > 0$ since the contact structure is Stein fillable. For $m = 0$, we have the empty monodromy. Stabilizing once, we see that this is Stein fillable, hence tight. For $m<0$, $c(\Sigma, \y^m)$ is overtwisted by Theorem \ref{thm:OvertwistedSoberingArc}.\\

\item[D:] If $d>0$, then $\by^d*\w*\y^m = \by^{d-1}*\w*\w*\y^2\x\y^2\x\y^m = \by^{d-1}*\w*\y^2\x\y^2*\x\y^2\x\y^2*\x\y^m$. Since $c(\Sigma, \x\y^2\x\y^2*\x\y^m) \neq 0$ for any $m$, we conclude by Theorem \ref{thm:NaturalityOfContactInvariant} that $c(\Sigma, \by^d*\y^m) \neq 0$, hence this open book is tight. \\

If $d < 0$, then $\by^d*\w*\y^m = \by^{d+1}*\x^{-1}\y^{-2}\x^{-1}\y^{-2}*\y^m$ consists of left handed Dehn twists around $x$ with Dehn twists around $y$, and is therefore overtwisted by Theorem \ref{thm:OvertwistedSoberingArc}. \\

If $d = 0$, then $\by^d*\w*\y^m = \y^2\x\y^2\x*y^m$. Note that $\y^2\x\y^2\x* \x^{-1}y^m = \y^2\x\y^{2+m}$, and $(\Sigma, \y^2\x\y^{2+m}) = (\Sigma, \x\y^{4+m})$, which is Stein fillable for $m \geq -4$. So, since $c(\Sigma, \y^2\x\y^2\x* \x^{-1}y^m) \neq 0$ for $m \geq -4$, we can conclude by Theorem \ref{thm:NaturalityOfContactInvariant} that $c(\Sigma, \y^2\x\y^2\x*\y^m) \neq 0$ for $m \geq -4$. For the proof of tightness when $m<-4$, see section \ref{sec:TypeD}.\\

\item[E:] If $d \geq 1$, then $\by^d*\x^m\y^{-1} = \by^{d-1}*\y\x^2\y\x^2\y\x^2\y\x^2*\x^m\y^{-1} = \by^{d-1}*\y\x^2\y\x^2\y\x^2\y\x^{2+m}\y^{-1}$ And, $(\Sigma,  \by^{d-1}*\y\x^2\y\x^2\y\x^2\y\x^{2+m}\y^{-1}) = (\Sigma,  \by^{d-1}*\x^2\y\x^2\y\x^2\y\x^{2+m}) = (\Sigma,  \by^{d-1}*\y\x^2\y\x^2\y\x^{4+m}),$ which is Stein fillable for $m \geq -3$.\\

If $d \leq 0$, then $\by^d*\x^m\y^{-1}$ consists of left handed Dehn twists around $x$ with Dehn twists around $y$ for $m \in \{-1,-2,-3\}$, and is therefore overtwisted by Theorem \ref{thm:OvertwistedSoberingArc}. \\

\item[F:] If $d \geq 0$, then $\by^d*\w*\x^m\y^{-1} = \by^{d}*\y\x^2\y\x^2*\x^m\y^{-1} = \by^{d}*\y\x^2\y\x^{2+m}\y^{-1}$ And, $(\Sigma,  \by^{d}*\y\x^2\y\x^{2+m}\y^{-1}) = (\Sigma,  \by^{d}*\x^2\y\x^{2+m}) = (\Sigma,  \y\x^{4+m}),$ which is Stein fillable for $m \geq -3$.\\

If $d < 0$, then $\by^d*\w*\x^m\y^{-1} = \by^{d+1}*\x^{-1}\y^{-2}\x^{-1}\y^{-2}*\x^m\y^{-1}$ consists of left handed Dehn twists around $x$ with Dehn twists around $y$ for $m \in \{-1,-2,-3\}$, and is therefore overtwisted by Theorem \ref{thm:OvertwistedSoberingArc}. 
\end{description}
 
 \qed

\begin{remark}
\label{remark:VeerDehn}
In light of Remark \ref{remark:Classification}, the proof of Theorem \ref{thm:TightGeneral} shows that, for periodic monodromy, tight is equivalent to Stein-fillable. Compare this with the results detailed in section \ref{sec:TightDehn}.
\end{remark}

\begin{remark}
\label{remark:LegendrianSurgery}
The proof also shows that tightness is equivalent to the non-vanishing of the contact invariant. In particular this implies that, for tight contact structures compatible with these open books, contact $-1$-surgery on a Legendrian knot is also tight. This follows from Theorem \ref{thm:NaturalityOfContactInvariant} and the fact that a Legendrian knot can be isotoped so that it lies in a page of the open book so that the contact framing is equal to the framing induced by the page.
\end{remark}

\begin{remark}It is also interesting to examine our results in the context of the following question, posed by Ozsv{\'a}th and Szab{\'o} in \cite{OSz1}: For a fibered knot $K \subset Y$, and $n$ large enough, it is clear that induced open book on $Y_{-1/n}(K)$ is Stein-fillable.  But what is the \emph{minimal} value of $n$ for which the induced open book on $Y_{-1/n}(K)$ is tight? In the genus one case, we have answered this question in Theorem \ref{thm:TightGeneral} for all monodromies.
\end{remark}

\newpage
\section{$Spin^c$ structures, Hopf invariants, and Type $D$ monodromies}
\label{sec:TypeD}
\begin{figure}[!htbp]
\begin{center}
\includegraphics[height=8cm]{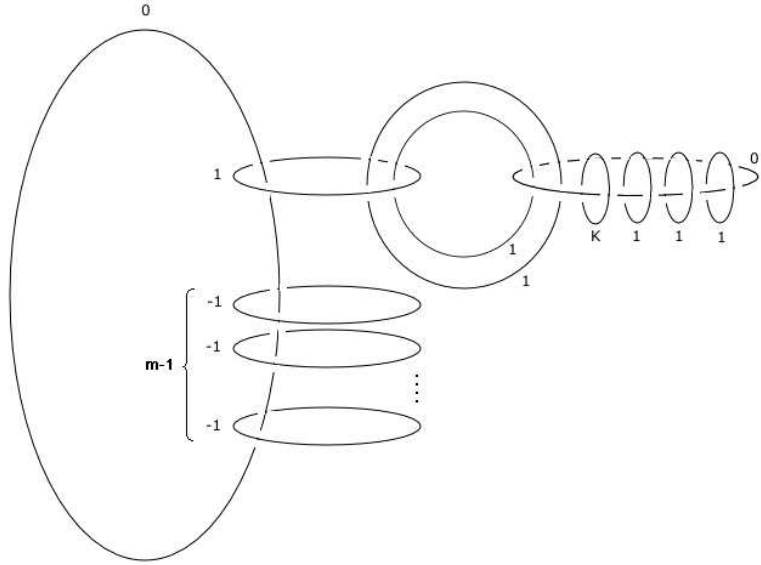}
\caption{\quad Surgery diagram for $-\widehat{M}(\w*\x\y^{-m})$. }
\label{fig:TypeD}
\end{center}
\end{figure}

In this section, we complete the proof of Theorem \ref{thm:TightGeneral} by showing that $c(\Sigma, \w*\y^{-m}) \neq 0$ for $m>4$. As we shall see, the characterization of tightness for type $D$ monodromies is a bit more involved than for the other types. If $\phi$ is a monodromy, let $M(\phi)$ denote the open book $(\Sigma,\phi)$, and let $\widehat{M}(\phi)$ denote the three-manifold underlying this open book. Figure \ref{fig:TypeD} is a surgery diagram for $$-\widehat{M}(\w*\x\y^{-m}) = -\widehat{M}(\y^{-m}*\w*\x) = -\widehat{M}(\y^{-m}*\y\x^2\y\x^2*\x)=-\widehat{M}(\y^{-(m-1)}*\x^2\y\x^2*\x). $$ 

\noindent Performing $-1$ and $0$ surgeries on the knot $K$, we obtain an exact triangle\\

\[
\begin{graph}(6,2)
\graphlinecolour{1}\grapharrowtype{2}
\textnode {A}(-1,1.5){$\heeg(-\widehat{M}(\w*\x\y^{-m}))$}
\textnode {B}(7, 1.5){$\heeg(-\widehat{M}(\w*\y^{-m})) $}
\textnode {C}(3, -.5){$\heeg(L(m,1))$}
\diredge {A}{B}[\graphlinecolour{0}]
\diredge {B}{C}[\graphlinecolour{0}]
\diredge {C}{A}[\graphlinecolour{0}]
\freetext (3,1.8){$\fwa$}
\freetext (0,.2){$\fwc$}
\freetext (6,.2){$\fwb$}
\end{graph}
\] \\

\noindent with the result that $\fwa$ maps $c(M(\w*\x\y^{-m})) \mapsto c(M(\w*\y^{-m}))$, by the naturality of the contact invariant (Theorem \ref{thm:NaturalityOfContactInvariant}). We have already proved that  $c(M(\w*\x\y^{-m})) \neq 0$ for all $m>0$ and that its image under $\fwa$ is non-trivial for $m \leq 4$ (since we have shown that $c(M(\w*\y^{-m})) \neq 0$ for $m \leq 4$). If we can show that its image under $\fwa$ is non-trivial for any $m$, then we will have completed the proof of Theorem \ref{thm:TightGeneral} for monodromies of type $D$. We prove the equivalent fact that $c(M(\w*\x\y^{-m})) \notin Im(\fwc)$. The proof proceeds as follows:

\begin{description}
\item[Step 1:] We show that the absolute grading of the contact invariant $c(M(\w*\x\y^{-m}))$ is $\frac{m-5}{4}$.
\item[Step 2:] We compute the grading shift of the maps $F_{W_3,s}$, where $s$ runs over all $spin^c$ structures on the cobordism $W_3$.
\item[Step 3:] Recall that the absolute gradings of the elements of $\widehat{HF}(L(m,1))$ are given by $\frac{(2j-m)^2-m}{4m}$ for $0 \leq j < m$. Combining this with the computations in steps $1$ and $2$, we show that $c(M(\w*\x\y^{-m}))$ cannot possibly be in the image of $F_{W_3}$ for grading reasons.
\end{description}

\subsection{Step 1}
\label{ssec:Step1}
\noindent Our first observation in the proof of Step 1 is the following:

\begin{thm}
\label{thm:SpinSpinc}
Let $\phi$ be the monodromy given by $\phi = \w*\x\y^{-m}$, and let $Y$ denote $\widehat{M}(\phi)$. For $m>0$, $c(M(\phi)) \in \heeg(-Y, s_c)$ where $s_c$ is a self-conjugate $spin^c$ structure on $-Y$. That is, $s_c = \bar{s_c}$.
\end{thm}

\begin{figure}[!htbp]
\begin{center}
\includegraphics[height=4cm]{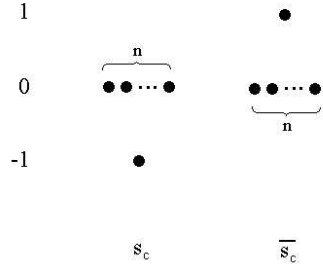}
\caption{\quad Possible Knot Floer homology in the $spin^c$ structures $s_c$ and $\bar{s_c}$. The filtration is given on the left. The dots represent generators.}
\label{fig:HFK}
\end{center}
\end{figure}

\noindent \textbf{Proof of Theorem \ref{thm:SpinSpinc}.}
This follows from examining the knot Floer homology $\widehat{HFK}(-Y, K)$, where $K$ is the binding of the open book decomposition $M(\phi)$ of $Y$. $\widehat{HFK}(-Y, K, -1)$ is generated by a single element which represents $c(M(\phi))$, and which is non-trivial in  $\heeg(-Y, s_c)$. There is a conjugation symmetry in knot Floer homology \cite{OSz3} which tells us that $\widehat{HFK}(-Y,K,1)$ is generated by a single element in the conjugate $spin^c$ structure $\bar{s_c}$. If $s_c \neq \bar{s_c}$, then the knot Floer homology in the $spin^c$ structures $s_c$ and $\bar{s_c}$ must look like the picture in Figure \ref{fig:HFK}.

\begin{figure}[!htbp]
\begin{center}
\includegraphics[height=4cm]{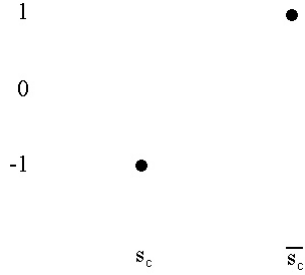}
\caption{\quad Another possibility for the knot Floer homology?}
\label{fig:HFK2}
\end{center}
\end{figure}

However, since knot Floer homology is the $E^1$ term in a spectral sequence that converges to $\heeg(-Y)$, and since $c(M(\phi))$ is the sole generator of $\heeg(-Y, s_c)$ (recall that $-\widehat{M}(\w*\x\y^{-m})$ is an $L$-space for $m>0$), there cannot be any non-trivial elements of $\widehat{HFK}(-Y,K,0)$ in the $spin^c$ structure $s_c$. By conjugation symmetry, this implies that there cannot be any non-trivial elements of $\widehat{HFK}(-Y,K,0)$ in the $spin^c$ structure $\bar{s_c}$. That is to say $n=0$ and the knot Floer homology in these two $spin^c$ structures must look like that depicted in Figure \ref{fig:HFK2}.

Yet, this last picture is not possible either, for there is also a symmetry under orientation reversal, $\widehat{HFK}_d(Y,K) = \widehat{HFK}^{-d}(-Y,K)$, which respects $spin^c$ structures, but reverses the sign of the filtration. If our knot Floer homology $\widehat{HFK}(-Y,K)$ looks like that depicted in Figure \ref{fig:HFK2}, then $\widehat{HFK}(Y,K)$ looks the same, only with the filtrations changing sign. This follows from the orientation reversal symmetry and Universal Coefficient Theorem, as we are using $\mathbb{Z}_2$ coefficients. Then the map $$\widehat{HFK}(Y,K,-1) \rightarrow \heeg(Y)$$ is non-trivial: there cannot be any higher differentials in the spectral sequence, so the generator of $\widehat{HFK}(Y,K,-1)$ must survive in $\heeg(Y)$. However, this is the statement that the contact invariant $c(M(\phi^{-1}))$ for the corresponding contact structure on $-Y$ is non-trivial. Yet, this contact structure is overtwisted by Goodman's sobering arc criterion since $\phi^{-1} = \y^m\x^{-1}*\y^{-2}\x^{-1}\y^{-2}\x^{-1}$, so we arrive at a contradiction.

\begin{figure}[!htbp]
\begin{center}
\includegraphics[height=4cm]{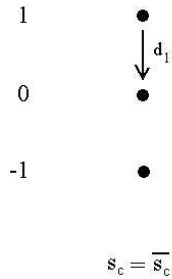}
\caption{\quad The knot Floer homology in the $spin^c$ structure $s_c = \bar{s_c}$}
\label{fig:HFK3}
\end{center}
\end{figure}

Consequently, it must be the case that $s_c = \bar{s_c}$. This finishes the proof of Theorem \ref{thm:SpinSpinc}. With a bit more work, it is possible to show that $\widehat{HFK}(-Y,K)$ looks like that depicted in Figure \ref{fig:HFK3} in the $spin^c$ structure $s_c = \bar{s_c}$. The $d_1$ differentials are shown in this figure and there can be no higher differentials in the related spectral sequence. 

\qed \\

To complete Step 1, we compute the absolute grading of $c(M(\w*\x\y^{-m}))$. Recall, from section \ref{sec:ContactInvariant}, that we have the following exact triangle
\[
\begin{graph}(6,2)
\graphlinecolour{1}\grapharrowtype{2}
\textnode {A}(-1,1.5){$\heeg(-\widehat{M}(\w*\x\y^{-m}))$}
\textnode {B}(7, 1.5){$\heeg(-\widehat{M}(\w*\x\y^{-m-1})) $}
\textnode {C}(3, -.5){$\heeg(Q(1)=S^3)$}
\diredge {A}{B}[\graphlinecolour{0}]
\diredge {B}{C}[\graphlinecolour{0}]
\diredge {C}{A}[\graphlinecolour{0}]
\freetext (3,1.8){$F_W$}
\end{graph}
\]
\\

\noindent where $F_W$ (summing over all $spin^c$ structures) is injective and maps $c(M(\w*\x\y^{-m})) \mapsto c(M(\w*\x\y^{-m-1}))$. Now, for a $spin^c$ structure $s$ on the cobordism $W$, $$F_{W,\bar{s}} = \mathfrak J F_{W,s} \mathfrak J$$ where $\mathfrak J:\heeg(X,t) \rightarrow \heeg(X, \overline{t})$ is the isomorphism on homology exhibited in section 3 of \cite{OSz5}. Therefore, $F_{W,\bar{s}}(c(M(\w*\x\y^{-m}))) =F_{W,s}(c(M(\w*\x\y^{-m})))$, as the $spin^c$ structure associated to $c(M(\w*\x\y^{-m}))$ is self-conjugate for any $m>0$. So, if $s \neq \bar{s}$, the contributions of these two maps cancel when we sum over $spin^c$ structures. In fact the only contributions which are not cancelled out are those coming from the maps $F_{W,s}$, where $s = \bar{s}$ on $W$. For such $s$, $c_1(s) = -c_1(s)$, hence $c_1(s)^2 = 0$. The grading shift of the map $F_{W,s}$ is given by $$\frac{c_1(s)^2 - 3\sigma -2\chi}{4}$$ which in this case is $1/4$: the cobordism is given by two-handle addition, so $\chi = 1$; in addition, this cobordism is negative definite (it cannot be positive definite as the map on $\heeg$ is non-zero, and it cannot be indefinite because all three terms in the associated surgery exact triangle are rational homology three-spheres), so $\sigma = -1$. Therefore, since $c(\w*\x\y^{-m}) \neq 0$ for $m>0$ and $c(M(\w*\x\y^{-m})) \mapsto c(M(\w*\x\y^{-m-1}))$, we may conclude that $$gr(c(M(\w*\x\y^{-m-1})))-gr(c(M(\w*\x\y^{-m}))) = 1/4.$$
For $m=1$, $-\widehat{M}(\w*\x\y^{-1}) = -L(5,1)$, and there is only one self-conjugate $spin^c$ structure $s_0$ on $-L(5,1)$. Moreover, the absolute grading of the generator of $\heeg(-L(5,1),s_0)$ is $-1$. Hence, $gr(c(M(\w*\x\y^{-1}))) = -1$, and we obtain the formula $gr(c(M(\w*\x\y^{-m}))) = \frac{m-5}{4}$ by induction, completing Step 1. By \cite{OSz1}, the absolute grading of the contact invariant is equal to the Hopf invariant of the corresponding two-plane field. Therefore, we have proved that the Hopf invariant of the two-plane field associated to the contact structure given by the monodromy $\phi = \w*\x\y^{-m}$ is $h(\phi) = \frac{m-5}{4}$. More generally, we can show that 

\begin{thm}
\label{thm:HopfInvariant}
The Hopf invariant of the two-plane field associated to the contact structure given by the monodromy $\phi = \w*\ls$ is $h(\phi) = -1+1/4*\sum_{i=1}^n{(b_i-a_i)}.$ Here, the $a_i,b_i \geq 0$ and $a_i \neq 0 \neq b_j$ for some $i,j$. 
\end{thm}

\noindent \textbf{Proof of Theorem \ref{thm:HopfInvariant}.}
The map $\heeg(-\widehat{M}(\w*\x^{a_1}\y^{-b_1}...\x^{a_n+1}\y^{-b_n})) \rightarrow \heeg(-\widehat{M}(\w*\ls))$ maps $c(M(\w*\x^{a_1}\y^{-b_1}...\x^{a_n+1}\y^{-b_n})) \mapsto c(M(\w*\ls)).$ For any monodromy of this form, the corresponding contact invariant lies is a self-conjugate $spin^c$ structure by the argument detailed above. The grading shift is therefore $1/4$, as before. Now Theorem \ref{thm:HopfInvariant} follows by induction. \qed

\subsection{Step 2}
\begin{figure}[!htbp]
\begin{center}
\includegraphics[height=7cm]{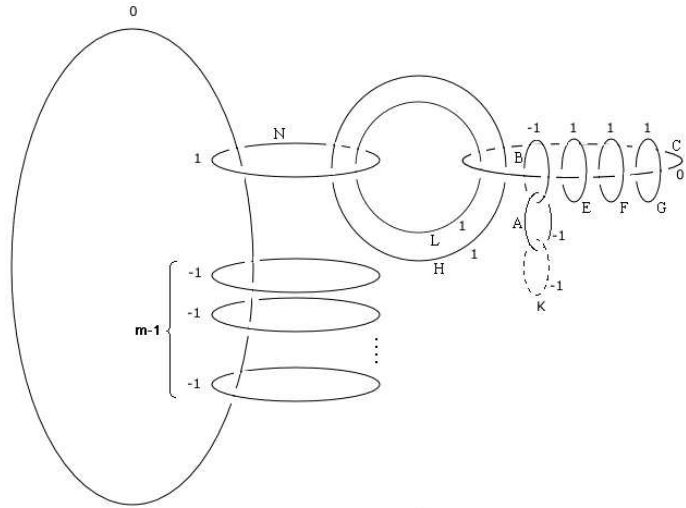}
\caption{\quad The cobordism $W_3$.}
\label{fig:Cobordism}
\end{center}
\end{figure}

The cobordism $W_3$ is given by attaching a $-1$ framed two-handle along $K$ in Figure \ref{fig:Cobordism}. The solid curves represent a cobordism, call it $Z$ from $S^3$ to $L(m,1)$, and the dashed curve $K$ represents the cobordism $W_3$. Now perform the following sequence of Kirby moves:

\begin{itemize}
\item{Blow down the $(m-1)$ $-1$-framed two-handles.}
\item{Blow down $A$, and slide $K$ over $C$.}
\item{Cancel $B$ against $C$}
\item{Blow down $E$, $F$, and $G$.}
\item{Blow down $H$, $L$, and $N$.}
\end{itemize}

\begin{figure}[!htbp]
\begin{center}
\includegraphics[height=3cm]{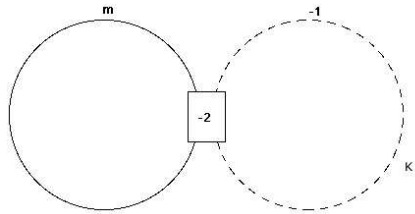}
\caption{\quad Another Kirby diagram for the cobordism $W_3$.}
\label{fig:Cobordism2}
\end{center}
\end{figure}

\noindent We are left with the Kirby Diagram in Figure \ref{fig:Cobordism2}. The solid curve represents a cobordism, call it $Z'$ from $S^3$ to $L(m,1)$. And the dashed curve represents $W_3$. Denote by $d(W,s)$ the grading shift of the map $F_{W,s}$ induced by a cobordism $W$ and a $spin^c$ structure $s$ on $W$. Again, $d(W,s) = \frac{c_1(s)^2 - 3\sigma(W) -2\chi(W)}{4}$. Under concatenation of cobordisms (where $spin^c$ structures agree on the common boundaries), this expression is additive. Denote the contatenation of $Z'$ with $W_3$ by $Z'*W_3$, which is a cobordism from $S^3$ to $-\widehat{M}(\w*\x\y^{-m})$. Then $d(W_3,s)$ = $d(Z'*W_3,t*s) - d(Z',t)$, where $t|_{L(m,1)} = s|_{L(m,1)}$. Alternatively, in order to compute $d(W_3,s)$, we simply compute $d(Z'*W_3, s') - d(Z', s'|_{Z'})$, where $s'|_{W_3} = s$. We take this approach and compute $d(Z'*W_3, s') - d(Z', s'|_{Z'})$ for all possible $s'$. From the formula for grading shift, we have $d(Z'*W_3, s') - d(Z', s'|{Z'})$ 
\begin{eqnarray*}
=&&\frac{[c_1(s')^2 - c_1(s'|_{Z'})^2] - 3[\sigma(Z'*W_3)-\sigma(Z')] -2[\chi(Z'*W_3)-\chi(Z')]}{4}\\
=&&\frac{c_1(s')^2 - c_1(s'|_{Z'})^2}{4} + \frac{1}{4}.
\end{eqnarray*} 

$Spin^c$ structures on a four-manifold $W$ are in one-to-one correspondence with characteristic vectors in $H^2(W;\mathbb{Z})$. By characteristic vector, we mean a cohomology class $K$ whose evaluation on a homology class $S \in H_2(W; \mathbb{Z})$ satisfies $<K,S> \equiv S^2\ mod\ 2$. This correspondence is given by $$s \mapsto c_1(s).$$ To understand $c_1(s)^2$, we need to think of it as a class in $H^2(W, \partial W; \mathbb{Q})$, where the intersection form on cohomology is defined. That is, we pull it back under the natural map $$H^2(W, \partial W; \mathbb{Q}) \rightarrow H^2(W; \mathbb{Q}).$$ and then compute the square of this pullback. Let $S_1,...,S_n$ be a basis for $H_2(W;\mathbb{Z})$. Then this map, when written as a matrix with respect to the bases $\{PD[S_1],...,PD[S_n]\}$ for $H^2(W,\partial W; \mathbb{Q})$ and $\{S_1^*,...,S_n^*\}$ for $H^2(W;\mathbb{Q})$, looks exactly like the intersection matrix for $H_2(W;\mathbb{Z})$.

The intersection matrix for the cobordism $Z'$ is is $(m)$, and the intersection matrix for $Z'*W_3$ is 
\[ \left( \begin{array}{cc}
m&-2\\
-2&-1
\end{array}\right)\]
Characteristic vectors in $H^2(Z'*W_3;\mathbb{Z})$ are of the form $(a,b)$, where $a \equiv m\ mod\ 2$ and $b \equiv 1\ mod\ 2$, when written with respect to the Hom-dual basis. It is clear that such a vector restricts to a characteristic vector in $H^2(Z';\mathbb{Z})$, and the restriction is simply $(a)$. Following the directions above, if $s$ on $Z'*W_3$ corresponds to the characteristic vector $(a,b)$, then $$c_1(s)^2 =\frac{a^2-4ab-mb^2}{m+4}$$ and $c_1(s|_{Z'})^2 = \frac{a^2}{m}$. And $$\frac{c_1(s')^2 - c_1(s'|_{Z'})^2}{4} = \frac{-m^2b^2-4abm-4a^2}{4m(m+4)}.$$ This completes Step 2.
\qed\\

\subsection{Step 3}
The absolute gradings for the various generators of $L(m,1)$ are given by $\frac{(2j-m)^2-m}{4m}$, where $j$ ranges from $0$ to $m-1$ \cite{OSz6}. Thus, if $F_{W_3}(x) = c(M(\w*\x\y^{-m}))$ for some $x \in \heeg(L(m,1))$ and some $m$, then $$\frac{(2j-m)^2-m}{4m} + \frac{-m^2b^2-4abm-4a^2}{4m(m+4)} = \frac{m-5}{4}$$ for some $0\leq j <m$, some $a$, and some $b$. But when $m$ is odd, this is impossible. Multiplying both sides by $4m$ and simplifying, we find that $$\frac{-m^2b^2-4abm-4a^2}{m+4} = -4m+4mj-4j^2.$$ However, the right hand side is divisible by $4$, hence $-m^2b^2-4abm-4a^2$ is also divisible by $4$. Then $m^2b^2$ is divisible by $4$. But recall that $b \equiv 1\ mod\ 2$, so if $m$ is odd, this cannot be. Thus, we have shown that $c(M(\w*\y^{-m}))$ is non-trivial when $m>0$ is odd. From this, we can draw the same conclusion for all $m>0$. Simply observe that the map $\heeg(-\widehat{M}(\w*\y^{-m})) \rightarrow \heeg(-\widehat{M}(\w*\y^{-m-1}))$ takes $c(M(\w*\y^{-m})) \mapsto c(M(\w*\y^{-m-1}))$. This completes the proof of Theorem \ref{thm:TightGeneral}. Whew!
\qed

\newpage
\section{$\Tight$ versus $\Dehn$}
\label{sec:TightDehn}
Let $\Tight$ denote the set of monodromies on $\Sigma$ which correspond to tight contact structures, and let $\Dehn$ denote the set of monodromies whose conjugate (by some element in $\mcg$) can be expressed as the product of right-handed Dehn twists. As was mentioned before, a contact structure is Stein-fillable if and only if it is compatible with \emph{some} open book whose monodromy is expressible as a product of right-handed Dehn twists. We should point out that there do exist Stein-fillable contact structures which are compatible with genus one, one boundary component open books, but whose monodromy cannot be taken to be in $\Dehn$. One such example is the contact structure compatible with the empty monodromy on $\Sigma$. After one stabilization, we see that this contact structure is Stein-fillable. On the other hand, it follows from the lemma below that the empty monodromy cannot be written as a product of right-handed Dehn twists around curves in $\Sigma$. Despite this discrepancy between Stein-fillable and $\Dehn$, an analysis of $\Tight - \Dehn$ seems to be an appropriate first step in the identification of tight but non-Stein-fillable contact structures which are compatible with genus one, one boundary component open books. 

This program was initiated by Honda, Kazez, and Mati{\'c} in \cite{HKM2}. In their paper, the authors do not explicitly study $\Tight$. Instead, they investigate $\Veer$, which is the monoid of right-veering diffeomorphisms of $\Sigma$. This is motivated by their recent results suggesting that $\Veer = \Tight$ when $\Sigma$ is the once-punctured torus (this has been verified in the case of pseudo-Anosov monodromies). The authors are able to identify infinitely many pseudo-Anosov monodromies in $\Veer - \Tight$. Their analysis makes use of the fact that a genus one, one boundary component open book is the branched double cover of a three-braid. Then, via a combination of the \emph{Rademacher function} and the \emph{rotation number}, they find three-braids which correspond to monodromies in $\Veer - \Dehn$.

Our result very closely parallels that of \cite{HKM2} but it is simpler in its statement and proof. We need only one lemma to identify infinitely many monodromies in $\Tight - \Dehn$:

\begin{lem}
\label{lem:Dehn}
If a monodromy $\phi = \x^{a_1}\y^{b_1}...\x^{a_n}\y^{b_n}$ is in $\Dehn$, then $\sum_{i=1}^n {a_i+b_i} > 0$. Moreover, if $\phi$ can be written as the product of $k$ right-handed Dehn twists about homologically non-trivial simple closed curves in $\Sigma$, then $\sum_{i=1}^n {a_i+b_i} =k$.
\end{lem}

Writing $\delta = (\x\y)^6$, we can easily identify infinitely many $\phi \in \Tight - \Dehn$ by combining Lemma \ref{lem:Dehn} with Theorem \ref{thm:TightGeneral}. For instance, monodromies of the form $\phi = \by^k * \w * \ls$ are in $\Tight$ when the $a_i,b_i \geq 0$, $k \geq 0$, and $a_i \neq 0 \neq b_j$ for some $i,j$, by Theorem \ref{thm:TightGeneral}. On the other hand, according to Lemma \ref{lem:Dehn}, $\phi \notin \Dehn$ if $$6+12k +\sum_{i=1}^n {a_i} \leq  \sum_{i=1}^n {b_i}.$$

\noindent \textbf{Proof of Lemma \ref{lem:Dehn}.}
The lemma follows from the well-known fact that $f D_{\gamma}f^{-1} =D_{f(\gamma)}$ (see, for example, \cite{OS}) where, in this notation, $D_{\gamma}$ represents a right-handed Dehn twist around the curve $\gamma$ and $f:\Sigma \rightarrow \Sigma$ is an orientation-preserving homeomorphism. Keeping with this notation for right-handed Dehn twists (until the end of this section), we need only to check that, when $\alpha$ is a homologically non-trivial curve, $D_{\alpha} = D_x^{a_1}D_y^{b_1}...D_x^{a_n}D_y^{b_n}$ for some collection of $a_i,b_j$ such that  $\sum_{i=1}^n {a_i+b_i} =1$. This is all we need to check because the relations in $\mcg$ between words in $D_x$ and $D_y$ are generated by the relation $D_xD_yD_x = D_yD_xD_y$, which preserves sums of exponents. 

In order to check this, observe that since $\alpha$ is homologically non-trivial, we can find an orientation-preserving homeomorphism $f$ so that $f(x)$ is isotopic to $\alpha$. Then, we have that $f^{-1} D_x f = D_{\alpha}$. Writing $f$ as a word in the Dehn twists $D_x$ and $D_y$, we are done.

\qed

\newpage
\section{$L$-spaces and Genus One Fibered Knots}
\label{sec:Lspaces}
Between Remark \ref{remark:OtherLspaces} and Lemma \ref{lem:wwLspace}, we have identified three types of $L$-spaces with genus one, one boundary component open book decompositions. These are $\ww$, $-\ww$, and $\wa$. 
By comparing surgery diagrams, it can be shown that $-\ww$ is equal to the three-manifold $\widehat{M}(-1;2,2,-a_1,...,-a_m)$ for some set of $a_{j} \in \mathbb{Z}^{\geq 0}$ and some $a_i \neq 0$. We summarize these statements and more in the following theorem:

\begin{thm}
\label{thm:AllLspaces}
For $b_{j} \in \mathbb{Z}^{\geq 0}$ and some $b_i \neq 0$, the following are $L$-spaces. Conversely, if $Y$ is an $L$-space with a genus one, one boundary component open book decomposition with pseudo-Anosov monodromy, then $Y$ takes one of the following forms:
\begin{description}
\item[1.]$\wa$
\item[2.]$\ww$
\item[3.]$\wx$
\end{description}
\end{thm}




\noindent \textbf{Proof of Theorem \ref{thm:AllLspaces}.}
These manifolds are all $L$-spaces by arguments identical to those in the proof of Lemma \ref{lem:wwLspace}.
For the converse, let $M_{\phi}$ denote the mapping torus of $\phi:\Sigma \rightarrow \Sigma$. Then $M_\phi$ is a three-manifold with torus boundary. Let $\widehat{M}_{\phi}(p/q)$ denote the $p/q$ Dehn filling of $M_{\phi}$ with respect to some framing. Roberts \cite{Ro} shows that if $\phi$ is pseudo-Anosov, then for 

\begin{itemize}
\item{$trace(\phi)>2$ and all but one Dehn filling, $\widehat{M}_{\phi}(p/q)$ has a co-orientable taut foliation, and for}
\item{$trace(\phi)<-2$, and $p/q \in (-\infty, 1)$ with respect to the framing in \cite{Ro}, $\widehat{M}_{\phi}(p/q)$ has a co-orientable taut foliation.}
\end{itemize}

\noindent On the other hand, $L$-spaces have no co-orientable taut foliations \cite{OSz2}. 

By Theorem \ref{thm:WordReduction}, the three-manifolds with genus one, one boundary component open book decompositions whose monodromy is pseudo-Anosov with $trace >2$ are of the form $\widehat{M}(k;-b_1,...,-b_n)$, and those with $trace < -2$ are of the form $\widehat{M}(k;2,2,-b_1,...,-b_n)$. As $k$ varies, these manifolds correspond to different Dehn fillings of mapping tori of the sort mentioned above. For $trace >2$, the $k = 0$ filling is an $L$-space, and therefore has no co-orientable taut foliation. Then \cite{Ro} tells us that $\widehat{M}(k;-b_1,...,-b_n)$ has a co-orientable taut foliation for $k \neq 0$, and is therefore not an $L$-space. 

Suppose that $\phi$ is the monodromy of the open book  $M(0;2,2,-b_1,...,-b_n)$. In the framing in \cite{Ro}, the longitude is the oriented boundary of a page of the open book and a meridian is chosen which intersects this longitude once. So, a priori, we know that $\wx = \widehat{M}_{\phi}(1/m)$ for some $m$, with respect to this framing. Since $\wx$ is an $L$-space, it must be true that $\wx$ is equal to either $ \widehat{M}_{\phi}(1/1)$ or $\widehat{M}_{\phi}(1/0)$ with respect to this framing, and $\ww$ is equal to the other. Then, $\widehat{M}(-k;2,2,-b_1,...,-b_n)$ is equal to $\widehat{M}_{\phi}(1/k)$ or $\widehat{M}_{\phi}(1/1-k)$. In either case, Roberts' results tell us that $\widehat{M}(-k;2,2,-b_1,...,-b_n)$ must have a co-orientable taut foliation for $k \neq 1$ or $0$, and is therefore not an $L$-space.
\qed \\

These $L$-spaces can be used to manufacture an infinite family of hyperbolic three-manifolds with no co-orientable taut foliations. The first such examples were found in \cite{RSS}, see also \cite{CD}. Let $\phi$ be the monodromy of the open book $M(0;2,2,-b_1,...,-b_n)$. Since $\ww$ is an $L$-space and $+1$-surgery on the binding of $\ww$ is an $L$-space (as it is equal to $\wx$), we can prove, using the surgery exact triangle and inductive arguments as before, that

\begin{thm}
\label{thm:HyperbolicLspaces}
$p/q$-surgery on the binding of $\ww$ is an $L$-space for $p/q \geq 1$. 
\end{thm}

\noindent According to Thurston \cite{Th1}, if $\phi$ is pseudo-Anosov, then $M_{\phi}$ is hyperbolic. In addition, Thurston's Hyperbolic Dehn Surgery Theorem \cite{Th2} guarantees that all but finitely many Dehn fillings of $M_{\phi}$ are hyperbolic as well. Thus, all but finitely many of the $L$-spaces in Theorem \ref{thm:HyperbolicLspaces} are hyperbolic. This family is much larger than the family of examples found in \cite{RSS}. The examples in \cite{RSS} can be expressed as $p/q$-surgery on the binding of manifolds of the form $\widehat{M}(0;2,2,0,0,...,0,-1)$ for $p/q \geq 1$ and $p$ odd. These manifolds are obtained as surgeries on the components of the Borromean rings. More precisely, $p/q$-surgery on the binding of $\widehat{M}(0;2,2,\overbrace{0,0,...,0}^{m},-1)$ is the manifold $B(p/q,1,m+5)$. As such, the first homology of the manifolds in \cite{RSS} is generated by at most two elements. On the other hand, infinitely many of the manifolds in Theorem \ref{thm:HyperbolicLspaces} have first homology generated by three elements. For instance, $p/q$-surgery on the binding of $\widehat{M}(0;2,2,\overbrace{0,0,...,0}^{m},-2)$ is the manifold $B(p/q,2,m+3)$. To be fair, in \cite{RSS} the authors prove that these manifolds have no taut foliations whatsoever, co-orientable or otherwise.

In \cite{KMOSz}, the authors exhibit an infinite family of monopole $L$-spaces (defined similarly in terms of monopole Floer homology) which are given by rational surgeries on the components of the Borromean rings for which the surgery coefficient on each component is $\geq 1$. A priori, it is not evident that our construction supplies any \emph{new} $L$-spaces which cannot be expressed as surgeries on the Borromean rings, although it seems very likely that this is the case. I suspect that this can be verified by comparing the \emph{graded} Heegaard Floer homologies of these various $L$-spaces.

\begin{remark}
\label{LSpaceConjecture}
 Ozsv{\'a}th and Szab{\'o} conjecture that if $Y$ is an integral homology three-sphere and $Y$ is an $L$-space, then $Y = n(\Sigma(2,3,5))\ \# \ m(-\Sigma(2,3,5))$ for some integers $n$ and $m$. Examining the first homologies of the open books corresponding to the monodromies of types A-F, we find that if an integer homology three-sphere $Y$ contains a genus one fibered knot, then $Y$ is the result of some $1/n$ surgery on one of the trefoils or on the figure eight. (For a break down of $|H_1(Y)|$ by monodromy, see \cite{B}.) If we additionally assume that $Y$ is an $L$-space then it follows from Moser (\cite{Mos}) that $Y$ can only be $S^3$, $+1$ surgery on the right-handed trefoil, or $-1$ surgery on the left-handed trefoil. The latter two surgeries produce the manifolds $\Sigma(2,3,5)$ and $-\Sigma(2,3,5)$, respectively, verifying the conjecture for three-manifolds which contain a genus one fibered knot.
\end{remark}

\newpage
\section{Appendix}
\label{sec:Appendix}
Here, we illustrate the proof of Lemma \ref{lem:H1Sum}. The proof is not especially revealing, but we include it for the sake of completeness. \\

\noindent \textbf{Proof of Lemma \ref{lem:H1Sum}.}
Below is a surgery diagram for $-\ww$ 
\begin{figure}[!htbp]
\begin{center}
\includegraphics[height=8cm]{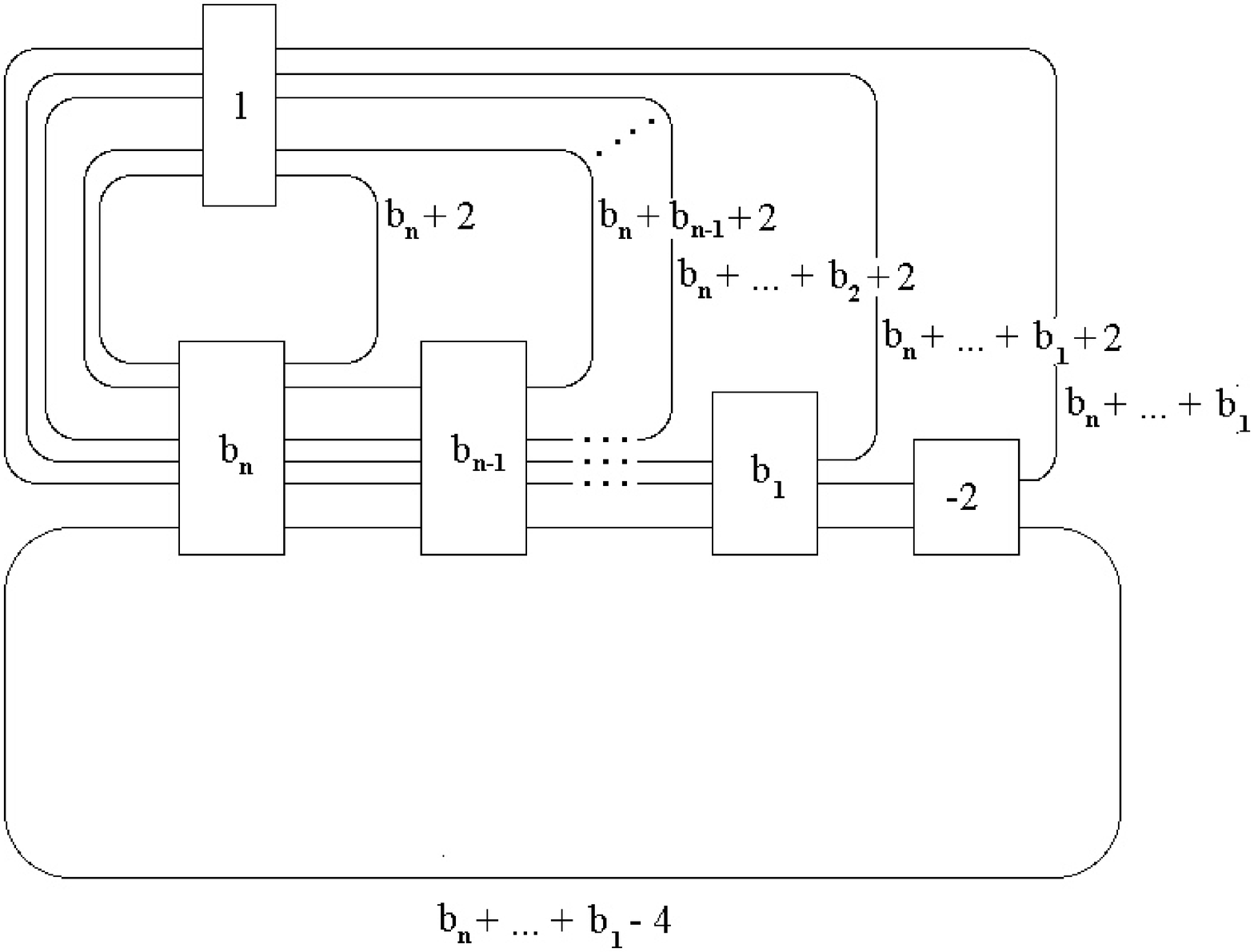}
\caption{\quad Surgery diagram for the 3--manifold $-\ww$}
\label{fig:BlownDownww}
\end{center}
\end{figure}

The linking matrix for $-\ww$ is given by 

\[ \left( \begin{array}{cccccc}
b_n+...+b_1-4 & b_n& b_n+b_{n-1}&  \cdots & b_n+...+b_1&  b_n+...+b_1-2 \\
b_n& b_n+2& b_n+1&\cdots&b_n+1&b_n+1\\
b_n+b_{n-1}& b_n+1&b_n+b_{n-1}+2&\cdots&b_n+b_{n-1}+1&b_n+b_{n-1}+1\\
\vdots & \vdots & \vdots & &\vdots & \vdots \\
b_n+...+b_1&b_n+1&b_n+b_{n-1}+1&\cdots&b_n+...+b_1+2&b_n+...+b_1+1\\
b_n+...+b_1-2&b_n+1&b_n+b_{n-1}+1&\cdots&b_n+...+b_1+1&b_n+...+b_1\\
\end{array} \right)\]

Then, the linking matrices for $-\wx$ and $\Qtwo$, respectively, are
\[ \left( \begin{array}{ccccccc}
b_n+...+b_1-4 & b_n& b_n+b_{n-1}&  \cdots & b_n+...+b_1&  b_n+...+b_1-2 &1\\
b_n& b_n+2& b_n+1&\cdots&b_n+1&b_n+1&1\\
b_n+b_{n-1}& b_n+1&b_n+b_{n-1}+2&\cdots&b_n+b_{n-1}+1&b_n+b_{n-1}+1&1\\
\vdots & \vdots & \vdots & &\vdots & \vdots \\
b_n+...+b_1&b_n+1&b_n+b_{n-1}+1&\cdots&b_n+...+b_1+2&b_n+...+b_1+1&1\\
b_n+...+b_1-2&b_n+1&b_n+b_{n-1}+1&\cdots&b_n+...+b_1+1&b_n+...+b_1&1\\
1&1&1&\cdots&1&1&-1\\
\end{array} \right)\]
and 
\[ \left( \begin{array}{ccccccc}
b_n+...+b_1-4 & b_n& b_n+b_{n-1}&  \cdots & b_n+...+b_1&  b_n+...+b_1-2 &1\\
b_n& b_n+2& b_n+1&\cdots&b_n+1&b_n+1&1\\
b_n+b_{n-1}& b_n+1&b_n+b_{n-1}+2&\cdots&b_n+b_{n-1}+1&b_n+b_{n-1}+1&1\\
\vdots & \vdots & \vdots & &\vdots & \vdots \\
b_n+...+b_1&b_n+1&b_n+b_{n-1}+1&\cdots&b_n+...+b_1+2&b_n+...+b_1+1&1\\
b_n+...+b_1-2&b_n+1&b_n+b_{n-1}+1&\cdots&b_n+...+b_1+1&b_n+...+b_1&1\\
1&1&1&\cdots&1&1&0\\
\end{array} \right)\]

\noindent Denote the determinants of these matrices by $\A$, $\An$, and $\Az$. It is clear that $$\An + \A = \Az.$$ 
\noindent By adding the last row to the previous rows, we can also see that $$\An = -A(b_1,...,b_n+1).$$
\noindent Observe that by performing row and column operations (begin by adding the last column to the first column and then the last row to the first row), we can write $\A = -C(b_1,...,b_n),$ where $C(b_1,...,b_n) =$
\[ \left| \begin{array}{cccccc}
b_n+2 & b_n+1 & b_n+1 & \cdots & b_n+1 & b_n+2\\
b_n+1 & b_n+b_{n-1}+2 & b_n+b_{n-1}+1 & \cdots & b_n+b_{n-1}+1 & b_n+b_{n-1}+2\\
b_n+1 & b_n+b_{n-1}+1 & b_n+b_{n-1}+b_{n-2}+2 & \cdots &  b_n+b_{n-1}+b_{n-2}+1 &  b_n+b_{n-1}+b_{n-2}+2\\
\vdots & \vdots & \vdots & & \vdots & \vdots\\
b_n+1 & b_n+b_{n-1}+1 & b_n+b_{n-1}+b_{n-2}+1 & \cdots & b_n+...+b_2+2 & b_n+...+b_2+2\\
b_n+2 & b_n+b_{n-1}+2 & b_n+b_{n-1}+b_{n-2}+2 & \cdots & b_n+...+b_2+2 & b_n+...+b_1+4\\
\end{array} \right|\]

\noindent We verify that $C(b_1,...,b_n) $ is positive and increasing in the parameter $b_n$, where the $b_{j} \in \mathbb{Z}^{\geq 0}$ and some $b_i \neq 0$. Once we establish this, it follows immediately that $$|\An| = |\A|+|\Az|,$$ which is equivalent to the statement of Lemma \ref{lem:H1Sum}. We start with two lemmas.

\begin{lem}
\label{lem:Positivity1}
For $b_1,...,b_n \geq 0$, $D(b_1,...,b_n) > 0$, where $D(b_1,...,b_n) =$
 \[ =\left| \begin{array}{cccccc}
b_n+1 & b_n & b_n & \cdots & b_n & b_n\\
b_n & b_n+b_{n-1}+1 & b_n+b_{n-1} & \cdots & b_n+b_{n-1} & b_n+b_{n-1}\\
b_n & b_n+b_{n-1} & b_n+b_{n-1}+b_{n-2}+1 & \cdots &  b_n+b_{n-1}+b_{n-2} &  b_n+b_{n-1}+b_{n-2}\\
\vdots & \vdots & \vdots & & \vdots & \vdots\\
b_n & b_n+b_{n-1} & b_n+b_{n-1}+b_{n-2} & \cdots & b_n+...+b_2 +1& b_n+...+b_2\\
b_n & b_n+b_{n-1} & b_n+b_{n-1}+b_{n-2} & \cdots & b_n+...+b_2 & b_n+...+b_1+1\\
\end{array} \right|\] 
\end{lem}

\noindent \textbf{Proof of Lemma \ref{lem:Positivity1}.}
This is clear in the case $n = 1$. We proceed by induction on $n$. $D(b_1,...,b_n) $
 \[ =\left| \begin{array}{ccccc}
b_n+b_{n-1}+1 & b_n+b_{n-1} & \cdots & b_n+b_{n-1} & b_n+b_{n-1}\\
b_n+b_{n-1} & b_n+b_{n-1}+b_{n-2}+1 & \cdots &  b_n+b_{n-1}+b_{n-2} &  b_n+b_{n-1}+b_{n-2}\\
\vdots & \vdots & & \vdots & \vdots\\
b_n+b_{n-1} & b_n+b_{n-1}+b_{n-2} & \cdots & b_n+...+b_2 +1& b_n+...+b_2\\
b_n+b_{n-1} & b_n+b_{n-1}+b_{n-2} & \cdots & b_n+...+b_2 & b_n+...+b_1+1\\
\end{array} \right|\]
 \[ +\left| \begin{array}{cccccc}
b_n & b_n & b_n & \cdots & b_n & b_n\\
b_n & b_n+b_{n-1}+1 & b_n+b_{n-1} & \cdots & b_n+b_{n-1} & b_n+b_{n-1}\\
b_n & b_n+b_{n-1} & b_n+b_{n-1}+b_{n-2}+1 & \cdots &  b_n+b_{n-1}+b_{n-2} &  b_n+b_{n-1}+b_{n-2}\\
\vdots & \vdots & \vdots & & \vdots & \vdots\\
b_n & b_n+b_{n-1} & b_n+b_{n-1}+b_{n-2} & \cdots & b_n+...+b_2 +1& b_n+...+b_2\\
b_n & b_n+b_{n-1} & b_n+b_{n-1}+b_{n-2} & \cdots & b_n+...+b_2 & b_n+...+b_1+1\\
\end{array} \right|\] 

\noindent The first summand is $ >0$ by induction. After row and column operations, we can write the second summand as 
 \[ =\left| \begin{array}{cccccc}
b_n & 0 & 0 & \cdots & 0 & 0\\
0 & b_{n-1}+1 & b_{n-1} & \cdots & b_{n-1} & b_{n-1}\\
0 & b_{n-1} & b_{n-1}+b_{n-2}+1 & \cdots &  b_{n-1}+b_{n-2} &  b_{n-1}+b_{n-2}\\
\vdots & \vdots & \vdots & & \vdots & \vdots\\
0 & b_{n-1} & b_{n-1}+b_{n-2} & \cdots & b_{n-1}+...+b_2 +1& b_{n-1}+...+b_2\\
0 & b_{n-1} & b_{n-1}+b_{n-2} & \cdots & b_{n-1}+...+b_2 & b_{n-1}+...+b_1+1\\
\end{array} \right|\] \\
\noindent And this is $\geq 0$ by induction. We have shown that we can write $$D(b_1,...,b_n) = D(b_1,..,b_{n-2},b_n+b_{n-1}) + b_n*D(b_1,...,b_{n-1}).$$ \noindent We use this fact in the next lemma.
\qed 

\begin{lem}
\label{lem:Positivity2}
For $b_1,...,b_n \geq 0$, $D(b_1,...,b_n) - D(b_2,...,b_n) \geq 0$.
\end{lem}
\noindent \textbf{Proof of Lemma \ref{lem:Positivity2}.}
For $n = 2$, this is clear. Again, we proceed by induction. We can write $D(b_1,...,b_n) - D(b_2,...,b_n)$
\begin{eqnarray*}
=  & D(b_1,..,b_{n-2},b_n+b_{n-1})  - D(b_2,..,b_{n-2},b_n+b_{n-1}) + b_n*[D(b_1,...,b_{n-1}) - D(b_2,...,b_{n-1})].
\end{eqnarray*} 
\noindent By induction, both $D(b_1,..,b_{n-2},b_n+b_{n-1})  - D(b_2,..,b_{n-2},b_n+b_{n-1})$ and $b_n*[D(b_1,...,b_{n-1}) - D(b_2,...,b_{n-1})]$ are $\geq 0$. \qed\\

Now we are ready to show that, for $b_1,...,b_n \geq 0$,  $C(b_1,...,b_n)$ is positive and increasing in the parameter $b_n$. This is clear for $n=1$. From here, we induct. $C(b_1,...,b_n) = $
\[ \left| \begin{array}{ccccc}
b_n+b_{n-1}+2 & b_n+b_{n-1}+1 & \cdots & b_n+b_{n-1}+1 & b_n+b_{n-1}+2\\
b_n+b_{n-1}+1 & b_n+b_{n-1}+b_{n-2}+2 & \cdots &  b_n+b_{n-1}+b_{n-2}+1 &  b_n+b_{n-1}+b_{n-2}+2\\
\vdots & \vdots & & \vdots & \vdots\\
b_n+b_{n-1}+1 & b_n+b_{n-1}+b_{n-2}+1 & \cdots & b_n+...+b_2+2 & b_n+...+b_2+2\\
b_n+b_{n-1}+2 & b_n+b_{n-1}+b_{n-2}+2 & \cdots & b_n+...+b_2+2 & b_n+...+b_1+4\\
\end{array} \right|\]
\[ +\left| \begin{array}{cccccc}
b_n+1 & b_n+1 & b_n+1 & \cdots & b_n+1 & b_n+2\\
b_n+1 & b_n+b_{n-1}+2 & b_n+b_{n-1}+1 & \cdots & b_n+b_{n-1}+1 & b_n+b_{n-1}+2\\
b_n+1 & b_n+b_{n-1}+1 & b_n+b_{n-1}+b_{n-2}+2 & \cdots &  b_n+b_{n-1}+b_{n-2}+1 &  b_n+b_{n-1}+b_{n-2}+2\\
\vdots & \vdots & \vdots & & \vdots & \vdots\\
b_n+1 & b_n+b_{n-1}+1 & b_n+b_{n-1}+b_{n-2}+1 & \cdots & b_n+...+b_2+2 & b_n+...+b_2+2\\
b_n+2 & b_n+b_{n-1}+2 & b_n+b_{n-1}+b_{n-2}+2 & \cdots & b_n+...+b_2+2 & b_n+...+b_1+4\\
\end{array} \right|\]

\noindent The first summand is $>0$ and increasing in the variable $b_n$ by induction. The second summand can be re-written as
\[ \left| \begin{array}{cccccc}
b_n+1 & 0 & 0 & \cdots & 0 & 1\\
0 & b_{n-1}+1 & b_{n-1} & \cdots & b_{n-1} & b_{n-1}\\
0 & b_{n-1} & b_{n-1}+b_{n-2}+1 & \cdots &  b_{n-1}+b_{n-2} &  b_{n-1}+b_{n-2}\\
\vdots & \vdots & \vdots & & \vdots & \vdots\\
0 & b_{n-1} & b_{n-1}+b_{n-2} & \cdots & b_n+...+b_2+1 & b_n+...+b_2\\
1 & b_{n-1} & b_{n-1}+b_{n-2} & \cdots & b_n+...+b_2 & b_n+...+b_1+1\\
\end{array} \right|\]
\noindent And this is equal to $(b_n+1)*D(b_1,...,b_{n-1}) - D(b_2,...,b_{n-1}) = b_n*D(b_1,...,b_{n-1}) + [D(b_1,...,b_{n-1}) - D(b_2,...,b_{n-1})]$. By Lemmas \ref{lem:Positivity1} and \ref{lem:Positivity2}, this is $\geq 0$ and increasing in $b_n$. This completes the proof of Lemma \ref{lem:H1Sum}. \qed

Let $\sigma_i$ be the $ith$ elementary symmetric polynomial in $n$ letters. If a bit more care is taken in our induction, it can be shown that $$|H_1(-\ww; \mathbb{Z})| = |A(b_1,...,b_n)| = 4 + \sum_{i=1}^n{(n-i+1)\sigma_i(b_1,...,b_n)}.$$

\end{document}